\documentclass[aop,MSNbibl,seceqn,dvips]{arximspdf}

\doi{10.1214/12-AOP754} 
\volume{41}
\issue{5}
\pubyear{2013}
\firstpage{3462}
\lastpage{3493}

\makeatletter
\newcommand{\rrVert}{\Vert}
\newcommand{\rrvert}{\vert}
\newcommand{\llVert}{\Vert}
\newcommand{\llvert}{\vert}

\newcommand{\Tilde}{\tilde}

\newcommand{\N}{\mathbb{N}}
\newcommand{\R}{\mathbb{R}}
\newcommand{\Z}{\mathbb{Z}}
\newcommand{\E}{\mathbb{E}}
\renewcommand{\P}{\mathbb{P}}
\renewcommand{\d}{{\mathrm{d}}}
\renewcommand{\L}{{\mathcal{L}_{t,\omega}}}

\newcommand{\var}{{\mathbb{V}\mathrm{ar}}_t^{\mu}}
\renewcommand{\O}{\mathcal{O}}

\newtheorem{theorem}{Theorem}[section]
\newtheorem{lemma}{Lemma}[section]
\newtheorem{proposition}{Proposition}[section]
\newtheorem{corollary}{Corollary}[section]

\newproclaim{remark}{Remark}[section]

\makeatother

\begin{document}
\begin{frontmatter}

\title{Annealed Brownian motion in a heavy tailed Poissonian potential}
\runtitle{Brownian motion in Poissonian potential}

\begin{aug}
\author[A]{\fnms{Ryoki} \snm{Fukushima}\corref{}\thanksref{t1}\ead[label=e1]{ryoki@math.titech.ac.jp}}
\runauthor{R. Fukushima}
\affiliation{Tokyo Institute of Technology}
\address[A]{Tokyo Institute of Technology\\
2-12-1, Ookayama\\
Meguro-ku, Tokyo 152-8551\\
Japan\\
\printead{e1}} 
\end{aug}

\thankstext{t1}{Supported by JSPS Grant-in-Aid for
Research Activity Start-up 22840019.}

\received{\smonth{10} \syear{2011}}
\revised{\smonth{3} \syear{2012}}

%
\begin{abstract}
Consider a $d$-dimensional Brownian motion in a random potential
defined by attaching a nonnegative and polynomially decaying potential
around Poisson points. We introduce a repulsive interaction between the
Brownian path and the Poisson points by weighting the measure by the
Feynman--Kac functional. We show that under the weighted measure, the
Brownian motion tends to localize around the origin. We also determine
the scaling limit of the path and also the limit shape of the random
potential.
\end{abstract}

%

\end{frontmatter}

\section{Introduction}\label{Intro}
\subsection{The model}
Consider a random potential defined by attaching the shape
function $\hat{v}(x)=|x|^{-\alpha}\wedge1$ ($d<\alpha<d+2$)
around a Poisson point process
$(\omega=\sum_i \delta_{\omega_i}, \P)$
with unit intensity as follows:
\[
V_{\omega}(x)=\sum_i \hat{v}(x-
\omega_i).
\]
Suppose we are also given the standard Brownian motion
$(\{X_t\}_{t \ge0}, \{P_x\}_{x\in\R^d})$ on~$\R^d$.
In this paper, we are interested in the long time behavior of a
Brownian motion under the \textit{annealed path measure}
defined by
\[
Q_t(A)=\frac{1}{Z_t} \E\otimes E_0 \biggl[ \exp
\biggl\{-\int_0^t V_{\omega}(X_s)
\,\d s \biggr\}\dvtx A \biggr],
\]
where $Z_t$ denotes the normalizing constant
\[
Z_t=\E\otimes E_0 \biggl[ \exp\biggl\{-\int
_0^t V_{\omega}(X_s) \,\d s \biggr\} \biggr].
\]
The weight $\exp\{-\int_0^t V_{\omega}(X_s) \,\d s\}$ introduces
a repulsive interaction between the Brownian path and Poisson points.
Since the averages are taken over both the path and the configuration,
it is natural to expect that $\omega$ tends to rarefy in a
region around the origin and the path favors to stay there.
However, it is often challenging to prove that such a localization
is typical under $Q_t$.

\subsection{Earlier studies}
We mention some early studies which are related to ours.
When $\alpha> d+2$, which is referred to as the light-tailed case,
Donsker and Varadhan \cite{DV75c} determined the asymptotics of
the normalizing constant
%
%
\begin{eqnarray}
\label{DV}
&&
\E\otimes E_0 \biggl[ \exp\biggl\{-\int
_0^t V_{\omega}(X_s) \,\d s \biggr\} \biggr]
\nonumber\\[-8pt]\\[-8pt]
&&\qquad = \exp\Bigl\{-\inf_{U\dvtx\mathrm{open}} \bigl\{\lambda_1(U)+|U|
\bigr
\} t^{{d}/({d+2})}\bigl(1+o(1)\bigr) \Bigr\}
\nonumber
\end{eqnarray}
as $t$ goes to $\infty$, where $|U|$ and $\lambda_1(U)$ stand for the
volume of $U$ and the smallest Dirichlet eigenvalue of
$-\Delta/2$ in $U$, respectively.
It follows from Faber--Krahn's inequality that
the unique minimizer of the above variational problem is
the ball with a certain radius $R_0$, up to translation.
This result
suggests that the dominant contribution to the right-hand side
of (\ref{DV}) comes from the following strategy:
there exists $x\in\R^d$ such that $\omega(B(x,R_0t^{1/(d+2)}))=0$ and
the Brownian motion $\{X_s\}_{0\le s\le t}$ stays in the ball.
Indeed, one can easily see that this specific event gives the
correct lower bound.
Motivated by this observation,
Sznitman \cite{Szn91b} ($d=2$) and Povel \cite{Pov99} ($d \ge3$)
proved that the above confinement is typical under the annealed
path measure when $\hat v$ has compact support. They also proved that the
scaled process $\{t^{-1/(d+2)}X_{t^{2/(d+2)}s}\}_{s \ge0}$
converges to a Brownian motion conditioned to stay in
a ball with radius $R_0$ and random center.
See also Bolthausen \cite{Bol94} for a
similar result in two-dimensional discrete space setting.

On the other hand, in the heavy tailed case $d < \alpha< d+2$,
Pastur \cite{Pas77} (first term) and Fukushima \cite{Fuk11a}
(second term) determined the asymptotics of the normalizing constant
%
%
\begin{eqnarray}
\label{Fuk10}
&&
\E\otimes E_0 \biggl[\exp\biggl\{-\int
_0^t V_{\omega}(X_s) \,\d s \biggr\} \biggr]\nonumber\\[-8pt]\\[-8pt]
&&\qquad = \exp\bigl\{-a_1t^{
{d}/{\alpha}} -
\bigl(a_2+o(1)\bigr)t^{({\alpha
+d-2})/({2\alpha})} \bigr\}\nonumber
\end{eqnarray}
as $t \to\infty$, where
%
%
\begin{eqnarray}
\label{var-Fuk}
a_1 &=& \omega_d \Gamma\biggl(
\frac{\alpha-d}{\alpha} \biggr),
\nonumber\\[-8pt]\\[-8pt]
a_2 &=& \mathop{\inf_{\phi\in W^{1,2}(\R^d),}}_{ \| \phi\|_{L^2}=1}
\biggl\{\int\frac{1}{2}\bigl|\nabla
\phi(x)\bigr|^2 + C(d,\alpha)|x|^2\phi(x)^2\,\d
x \biggr\}
\nonumber
\end{eqnarray}
with $\omega_d$ being the volume of the unit ball, $C(d,\alpha)>0$ a
constant and $W^{1,2}(\R^d)$ the usual Sobolev space. The constant
$C(d,\alpha)$ admits an explicit expression
$\frac{\alpha\sigma_d}{2d}\Gamma( \frac{2\alpha-d+2}{\alpha } )$ with
$\sigma_d$ the surface area of the unit sphere. As in the light-tailed
case, the correct lower bound of (\ref{Fuk10}) can
be given by considering a specific strategy:
\begin{longlist}[(St4)]
\item[(St1)] $V_{\omega}(0) = a_1\frac{d}{\alpha}t^{-({\alpha-d})/{\alpha}}
+o(t^{-({\alpha-d+2})/({2\alpha})})$;\vspace*{1pt}
\item[(St2)] $V_{\omega}(x)-V_{\omega}(0)= C(d,\alpha)
t^{-({\alpha-d+2})/{\alpha}} |x|^2
+o(t^{-({\alpha-d+2})/({2\alpha})})$
for\vspace*{1pt} $|x|<Mt^{({\alpha-d+2})/({4\alpha})}$;
\item[(St3)] $\sup_{0 \le s \le
t}|X_s|<Mt^{({\alpha-d+2})/({4\alpha})}$;\vspace*{1pt}
\item[(St4)] $L_t=\frac{1}{t}\int_0^t\delta_{X_s}\,\d s$
is, after the diffusive scaling with spatial factor
$t^{({\alpha-d+2})/({4\alpha})}$, weakly close to
$\phi_1(x)^2\,\d x$, where
\[
\phi_1(x)= \biggl(\frac{\sqrt{2C(d,\alpha)}}{\pi} \biggr)^{{d}/{4}}
\exp
\Biggl\{-\sqrt{\frac{C(d,\alpha)}{2}}|x|^2 \Biggr\}
\]
is the unique minimizer of (\ref{var-Fuk}).
\end{longlist}
Let us informally explain how this strategy
gives the correct lower bound.
The first event (St1) has probability
\begin{eqnarray*}
&&
\P\biggl(V_{\omega}(0) = a_1\frac{d}{\alpha} t^{-({\alpha-d})/\alpha}
+o \bigl(t^{-({\alpha-d+2})/({2\alpha})} \bigr) \biggr)\\
&&\qquad
\ge\exp\biggl\{-a_1
\frac{\alpha-d}{\alpha}t^{d/\alpha} +o \bigl(t^{({\alpha+d-2})/({2\alpha})} \bigr) \biggr\}.
\end{eqnarray*}
The second event (St2) introduces no extra cost since conditioned on
(St1), its probability is not too small.
In fact, it is close to 1 when $\alpha>2$, and
even when $\alpha\le2$, it is only polynomially small.
As for the third and fourth events,
Donsker and Varadhan's large deviation principle shows that
\begin{eqnarray}\label{LDP-lower}
&&
P_0\bigl(\mbox{(St3) and (St4) hold}\bigr) \\
&&\qquad
\ge\exp\biggl
\{-t^{({\alpha+d-2})/({2\alpha})} \int\frac
{1}{2}\bigl|\nabla\phi_1
(x)\bigr|^2\,\d x +o \bigl(t^{({\alpha+d-2})/({2\alpha})} \bigr)
\biggr\}\nonumber
\end{eqnarray}
after letting $M\to\infty$.
Summarizing the above, we obtain the correct lower bound
\begin{eqnarray*}
&&\E\otimes E_0 \biggl[\exp\biggl\{-\int_0^{t}
V_{\omega}(X_s)\,\d s \biggr\}\dvtx\mbox{(St1)--(St4) hold}
\biggr]
\\
&&\qquad \ge\exp\biggl\{-a_1t^{{d}/{\alpha}}\\
&&\qquad\quad\hspace*{19.8pt}{}
-t^{({\alpha+d-2})/({2\alpha})} \biggl(\int
\frac{1}{2}\bigl|\nabla\phi_1 (x)\bigr|^2 +C(d,
\alpha)|x|^2 \phi_1 (x)^2\,\d x+o(1) \biggr)
\biggr\}.
\end{eqnarray*}
From this observation, it is natural to ask whether the above
strategy is typical under $Q_t$ or not.
%
%
\begin{remark}\label{1.1}
Our model can be viewed as an example of diffusion in a spatially
correlated random potential. There are several results, mainly
concerning the asymptotics of $Z_t$, for such models.
See, for example, \cite{GM00} (discrete space) and \mbox{\cite{GK00,GKM00}}
(continuous space).
\end{remark}
%
\subsection{Results}
The results of the present paper show that the events
(St1)--(St4) consisting the optimal strategy in the last
subsection are, with appropriate modifications, typical
under the annealed path measure $Q_t$.
Moreover, the scaling limit of the path is also identified.
The first result is about the path localization (St3).
%
%
\begin{theorem}
\label{localization}
For any $\delta>0$,
\[
\lim_{t\to\infty} Q_t \Bigl(\sup_{0\le s\le t}|X_s|
\le t^{({\alpha-d+2})/({4\alpha})} (\log t)^{{1}/{2}+\delta} \Bigr)=1.
\]
\end{theorem}
%
%
%
\begin{remark}
The logarithmic correction corresponds to the
$M\to\infty$ operation in (\ref{LDP-lower}).
The power $1/2$ is natural in view of Theorem~\ref{scaling-limit} below.
\end{remark}
%
The second result says that the random potential viewed from its
local minimum looks like the quadratic function
\[
p_t(x)=C(d,\alpha)t^{-({\alpha-d+2})/{\alpha}}|x|^2,
\]
that is, (St2) is typical.
It also turns out that the barycenter $m_{L_t}$ of the
occupation time measure $L_t=\frac{1}{t}\int_0^t\delta_{X_s}\,\d s$
is close to a \textit{local minimizer} of $V_{\omega}$.
%
%
\begin{theorem}
\label{potential-confinement}
Let $m_t(\omega)$ be the point where $V_{\omega}$
attains the minimum in $B(0,t^{({\alpha-d+2})/({4\alpha})}\log t)$.
Then there exists $\varepsilon_0>0$ such that
%
%
\begin{eqnarray}
\label{thm2}
\hspace*{14pt}&&
\lim_{t\to\infty} Q_t \bigl(\sup\bigl
\{\bigl|V_{\omega}(x)-V_{\omega}\bigl(m_t(\omega)\bigr)
-p_t\bigl(x-m_t(\omega)\bigr)\bigr|\dvtx
\nonumber\\
&&\hspace*{93.7pt}\bigl|x-m_t(\omega)\bigr|\le t^{({\alpha-d+2})/({4\alpha})} \log t \bigr\}\le
t^{-({\alpha-d+2})/({2\alpha})-\varepsilon_0} \bigr)\\
&&\qquad= 1.
\nonumber
\end{eqnarray}
Moreover,
%
%
\begin{equation}
\label{bottom-location} \lim_{t\to\infty} Q_t
\bigl(\bigl|m_t(\omega)-m_{L_t}\bigr|<\varepsilon t^{({\alpha-d+2})/({4\alpha})}
\bigr)=1
\end{equation}
for any $\varepsilon>0$.
\end{theorem}
%
%
%
\begin{remark}
A statement similar to (\ref{thm2}) is proved by Gr\"{u}ninger and
K\"{o}nig \cite{GK09} for the parabolic Anderson model on $\Z^d$
with so-called ``almost bounded potentials.''
They call this phenomenon ``potential confinement.''\vadjust{\goodbreak}
\end{remark}
From Theorems~\ref{localization} and~\ref{potential-confinement},
we can deduce
scaling limits of the path, occupation time measure,
and local minimum of the potential.
We introduce the scale function $r(t)=t^{({\alpha-d+2})/({4\alpha})}$
and define the scaled process by $\Tilde{X}_s=r(t)^{-1}X_{r(t)^2s}$.
(The dependence of $\Tilde{X}$ on $t$ is omitted,
but this does not seem to cause any confusion.)
Let us begin with the result for the occupation time measure
$\Tilde{L}_t=\frac{1}{tr(t)^{-2}}\int_0^{tr(t)^{-2}}\delta_{\Tilde
{X}_s}\,\d s$,
which implies that (St4) is typical.
Let $\Tilde{m}_t(\omega)=r(t)^{-1}m_t(\omega)$ and
$\nu_m$ denote the measure with density
\[
\phi_1(x-m)^2 = \biggl(\frac{\sqrt{2C(d,\alpha)}}{\pi}
\biggr)^{{d}/{2}} \exp\bigl\{-\sqrt{2C(d,\alpha)}|x-m|^2 \bigr
\}.
\]

%
%
\begin{theorem}
\label{occupation-time}
For any $f \in C_c(\R^d)$ and $\varepsilon>0$,
%
%
\begin{equation}
\label{thm4-1} \lim_{t \to\infty}Q_t \biggl(\biggl\llvert\int f \,\d
\Tilde{L}_t -\int f \,\d\nu_{\Tilde{m}_t(\omega)}\biggr\rrvert
<\varepsilon
\biggr)=1.
\end{equation}
\end{theorem}
%
Next we state results on the local minimum of $V_{\omega}$.
It turns out that there is a gap between the expected value of
$V_{\omega}(m_t(\omega))$ and the right-hand side of (St1).
Moreover, the fluctuation of $V_{\omega}(m_t(\omega))$
is even larger than the gap when $\alpha< 2$.
%
%
\begin{theorem}
\label{local-minimum}
The following hold:
\begin{longlist}
\item$Q_t[V_{\omega}(m_t(\omega))]
=a_1\frac{d}{\alpha}t^{-({\alpha-d})/{\alpha}}
- (\sqrt{\frac{C(d,\alpha)}{8}}+o(1) )
t^{-({\alpha-d+2})/({2\alpha})}$;
\item$t^{({2\alpha-d})/({2\alpha})}
\{V_{\omega}(m_t(\omega))-Q_t[V_{\omega}(m_t(\omega))]\}$
converges\vspace*{1pt} in law to the
Gaussian random variable with variance
$\alpha\sigma_d\Gamma(\frac{3\alpha-d+1}{\alpha})$.
\end{longlist}
\end{theorem}
We finally state the result on the scaling limit of the path.
We write $R^m_0$ for the law of the Ornstein--Uhlenbeck process
starting at the origin with generator
$-\Delta/2+\sqrt{2C(d,\alpha)}\langle x-m,\nabla\rangle$.
We call $m$ the \textit{center}.
%
%
\begin{theorem}
\label{scaling-limit}
The process $\{\Tilde{X}_s\}_{s \ge0}$ under $Q_t$
converges as $t\to\infty$ in law to the Ornstein--Uhlenbeck
processes with random centers
\[
\int\d m \biggl(\frac{\sqrt{C(d,\alpha)}}{\sqrt{2}\pi} \biggr)^{
{d}/{2}} \exp\Biggl\{-\sqrt{
\frac{C(d,\alpha)}{2}}|m|^2 \Biggr\} R^m_0.
\]
\end{theorem}
%
%
%
\begin{remark}
It becomes clear in the proof that the random center $m$
corresponds to $\Tilde{m}_t(\omega)=r(t)^{-1}m_t(\omega)$.
Hence roughly speaking, this theorem means that
$\{\Tilde{X}_s\}_{s \ge0}$ behaves
like the Ornstein--Uhlenbeck process centered at
$\Tilde{m}_t(\omega)$.
\end{remark}
%
\subsection{Organization of the paper}
We briefly comment on the outline of the proof and
the organization of the paper.
The main difficulty lies in the proof of Theorem~\ref{localization}
and Theorem~\ref{potential-confinement}. They are closely related\vadjust{\goodbreak}
in the sense that the path localization
implies the potential confinement and vice versa.
But of course we need another input to leave the circular argument.
The key idea is that we actually need only a weaker result than
Theorem~\ref{localization} to deduce the potential confinement.
We prove such a weak localization in
Section~\ref{sect-weak-localization} by using crude estimates.
Then in Section~\ref{weak-to-PC}, we show that it indeed implies
a potential confinement (Proposition~\ref{PC-L_t}),
which is a slight modification of
Theorem~\ref{potential-confinement}.
Given the potential confinement,
Theorem~\ref{localization} follows by repeating a part of argument
in Section~\ref{sect-weak-localization}, and it is done in
Section~\ref{strong-localization}.
In the end of Section~\ref{strong-localization},
we show (\ref{bottom-location}) of
Theorem~\ref{potential-confinement} by
combining Theorem~\ref{localization} and Proposition~\ref{PC-L_t},
and it completes the proof of Theorem~\ref{potential-confinement}.
Sections~\ref{OT},~\ref{fluctuation} and~\ref{SL} are devoted to
prove Theorems~\ref{occupation-time},~\ref{local-minimum}
and~\ref{scaling-limit}, respectively.\looseness=1
\subsection{Notation}
For functions $f, g\dvtx[0,\infty) \to\R$, we say
$f(t)=\O(\exp\{g(t) \})$ if there exists $c>0$ such that
$f(t) = O(\exp\{cg(t) \})$ as $t\to\infty$.
\section{Weak localization}
\label{sect-weak-localization}
In this section, we prove the following weaker version of
Theorem~\ref{localization}.
%
%
\begin{proposition}
\label{weak-localization}
There exists $M_1>0$ such that
\[
\lim_{t\to\infty}Q_t \Bigl(\sup_{0\le s \le t}
|X_s| \le M_1 t^{({\alpha-d+6})/({8\alpha})} \Bigr)=1.
\]
\end{proposition}
%
We use the notation
\[
h_t=a_1\frac{d}{\alpha}t^{-({\alpha-d})/{\alpha}},
\]
which is expected to be the typical height of the bottom of
$V_{\omega}$ under $Q_t$ as explained in the \hyperref[Intro]{Introduction}.
We first show that $\min_{0 \le s \le t} V_{\omega}(X_s)$
cannot be too far from $h_t$.
%
%
\begin{lemma}
\label{lem1}
For any $M\in\R$ and sufficiently large $t$,
\begin{eqnarray*}
&&\P\Bigl(\min_{x \in(-t,t)^d} V_{\omega}(x)-h_t <
Mt^{-({3\alpha-3d+2})/({4\alpha})} \Bigr)
\\
&&\qquad \le\exp\biggl\{-a_1\frac{\alpha-d}{\alpha}t^{d/\alpha
}
+Mt^{({\alpha+3d-2})/({4\alpha})} -|M|t^{({\alpha+d-2})/({2\alpha})} \biggr\}.
\end{eqnarray*}
\end{lemma}
%
%
\begin{pf}
The event in the left-hand side concerns infinitely many points,
but we can reduce it to finite points since $V_{\omega}$
is smooth where it is small.
Indeed, if $V_{\omega}(x)<1$, then obviously
\[
V_{\omega}(x)=\sum_i |x-
\omega_i|^{-\alpha}\wedge1 =\sum_i
|x-\omega_i|^{-\alpha},
\]
and in particular $\min_i |x-\omega_i|>1$. Since there exists $c_0>0$
such that\break $|\nabla|x|^{-\alpha}|<c_0|x|^{-\alpha}$ for $|x|>1$,
it follows that
\[
\bigl|\nabla V_{\omega}(x)\bigr| \le\sum_i \bigl|
\nabla|x-\omega_i|^{-\alpha}\bigr| \le\sum
_i c_3|x-\omega_i|^{-\alpha}
=c_0 V_{\omega}(x)<c_0,
\]
when $V_{\omega}(x)<1$.
Therefore one can find $N>0$ such that
\begin{eqnarray*}
&&\P\Bigl(\min_{x \in(-t,t)^d} V_{\omega}(x)-h_t <
Mt^{-({3\alpha-3d+2})/({4\alpha})} \Bigr)
\\
&&\qquad \le\P\Bigl(\min_{x \in(-t,t)^d \cap t^{-N}\Z^d} V_{\omega}(x)-h_t <
\bigl(M+t^{-1}\bigr)t^{-({3\alpha-3d+2})/({4\alpha})} \Bigr)
\\
&&\qquad \le(2t+1)^{dN}\P\bigl(V_{\omega}(0)-h_t <
\bigl(M+t^{-1}\bigr)t^{-({3\alpha-3d+2})/({4\alpha})} \bigr).
\end{eqnarray*}
To bound the probability in the last line, we use the
asymptotics of the moment generating function
%
%
\begin{equation}
\label{MGF} \E\bigl[\exp\bigl\{-sV_{\omega}(0)\bigr\}\bigr] =\exp\bigl
\{-a_1s^{{d}/{\alpha}}+O\bigl(e^{-s}\bigr) \bigr\}
\qquad\mbox{as }s \to\infty,
\end{equation}
which is proved in Lemma 1 of \cite{Fuk11a}.
Taking $s=t(1+t^{-({d+2-\alpha})/({4\alpha})})$ in (\ref{MGF})
and using Chebyshev's inequality, one can see by a straightforward
calculation that
\begin{eqnarray*}
&&\P\bigl(V_{\omega}(0)-h_t < \bigl(M+t^{-1}
\bigr)t^{-({3\alpha-3d+2})/({4\alpha})} \bigr)
\\
&&\qquad \le\exp\biggl\{-a_1\frac{\alpha-d}{\alpha}t^{d/\alpha
}
+Mt^{({\alpha+3d-2})/({4\alpha})} -\bigl(|M|+\varepsilon\bigr)t^{({\alpha+d-2})/({2\alpha})} \biggr\}
\end{eqnarray*}
for some small $\varepsilon>0$.
\end{pf}
%
%
%
\begin{lemma}
\label{lem2}
\begin{eqnarray*}
&&
Q_t \Bigl(\min_{0\le s \le t} V_{\omega}(X_s)-h_t
\notin t^{-({3\alpha-3d+2})/({4\alpha})} [-2a_2, 2a_2 ] \Bigr)\\
&&\qquad =\O\bigl(
\exp\bigl\{-t^{({\alpha+d-2})/({2\alpha})} \bigr\} \bigr).
\end{eqnarray*}
\end{lemma}
%
%
\begin{pf}
We shall only give the proof of
\[
Q_t \Bigl(\min_{0\le s \le t} V_{\omega}(X_s)-h_t
> 2a_2t^{-({3\alpha-3d+2})/({4\alpha})} \Bigr) =\O\bigl(\exp\bigl\{
-t^{({\alpha+d-2})/({2\alpha})}
\bigr\} \bigr).
\]
The other half can be shown by a similar argument.
(See also Lemma~\ref{lem6} where a stronger statement is proved.)
Note first that
\begin{equation}
\label{macro-box}
P_0 \bigl(X_{[0,t]}\not\subset(-t,t)^d
\bigr) =\O\bigl(\exp\{-t\} \bigr)
\end{equation}
by a simple application of the reflection principle and that
\begin{eqnarray*}
&&\E\otimes E_0 \biggl[\exp\biggl\{-\int_0^{t}
V_{\omega}(X_s)\,\d s \biggr\}\dvtx\min_{0\le s \le t}
V_{\omega}(X_s) >2a_1t^{-({\alpha-d})/{\alpha}} \biggr]
\\
&&\qquad \le\exp\bigl\{-2a_1t^{{d}/{\alpha}} \bigr\}.
\end{eqnarray*}
From these bound and (\ref{Fuk10}), it follows that
\begin{eqnarray*}
&&Q_t \Bigl(\bigl\{X_{[0,t]}\not\subset(-t,t)^d
\bigr\} \cup\Bigl\{\min_{0\le s \le t} V_{\omega}(X_s)
>2a_1t^{-({\alpha-d})/{\alpha}} \Bigr\} \Bigr)
\\
&&\qquad =\O\bigl(\exp\bigl\{-t^{({\alpha+d-2})/({2\alpha})} \bigr\} \bigr).
\end{eqnarray*}

Next, for each $k \in\N$, let us introduce two events
\begin{eqnarray*}
E_k &=& \Bigl\{\min_{0\le s \le t} V_{\omega}(X_s)
-h_t-2a_2t^{-({3\alpha-3d+2})/({4\alpha})} \in t^{-1}[k,k+1)
\Bigr\},
\\
F_k &=& \Bigl\{\min_{x \in(-t,t)^d} V_{\omega}(x)
-h_t<2a_2t^{-({3\alpha-3d+2})/({4\alpha})}+t^{-1}(k+1) \Bigr\}.
\end{eqnarray*}
Then, it follows by the fact $E_k\cap\{X_{[0,t]}\subset(-t,t)^d\}
\subset F_k$
and Lemma~\ref{lem1} that
\begin{eqnarray*}
&&\E\otimes E_0 \biggl[\exp\biggl\{-\int_0^{t}
V_{\omega}(X_s)\,\d s \biggr\}\dvtx E_k\cap\bigl
\{X_{[0,t]}\subset(-t,t)^d\bigr\} \biggr]
\\
&&\qquad \le\exp\biggl\{-a_1\frac{d}{\alpha}t^{{d}/{\alpha}}
-2a_2t^{({\alpha+3d-2})/({4\alpha})}-k \biggr\} \P(F_k)
\\
&&\qquad \le\exp\bigl\{-a_1t^{{d}/{\alpha}} -2a_2t^{({\alpha+d-2})/({2\alpha})}+1
\bigr\}.
\end{eqnarray*}
Dividing both sides by $Z_t$ and
summing over $k \in[0, 2a_1t^{d/\alpha}+1]\cap\Z$,
we find
\begin{eqnarray*}
&&
Q_t \Bigl(\bigl\{X_{[0,t]}\subset(-t,t)^d\bigr
\} \\
&&\quad\hspace*{5.6pt}{}\cap\Bigl\{\min_{0\le s \le t} V_{\omega}(X_s)-h_t
\in\bigl[2a_2t^{-({3\alpha-3d+2})/({4\alpha})}, 2a_1t^{-({\alpha-d})/{\alpha}}
\bigr] \Bigr\} \Bigr)
\\
&&\hspace*{5.6pt}\qquad =\O\bigl(\exp\bigl\{-t^{({\alpha+d-2})/({2\alpha})} \bigr\} \bigr),
\end{eqnarray*}
and the proof is complete.
\end{pf}
%
Let us define a \textit{low level set} of $V_{\omega}$ by
\[
\L= \bigl\{x \in\R^d\dvtx V_{\omega}(x) \le h_t+2a_2t^{-({3\alpha-3d+2})/({4\alpha})}
\bigr\}.
\]
Then Lemma~\ref{lem2} above shows that
$Q_t(X_{[0,t]}\cap\L\neq\varnothing) \to1$.
In what follows, we are going to show that if $x \in\L$,
then $V_{\omega}$ is bounded below by a certain quadratic
function in an annulus around $x$. This has two consequences
which lead us to the weak localization bound:
\begin{longlist}
\item each connected component of $\L$ is not too large;
\item once the Brownian motion hits a connected component of
$\L$, it is difficult to escape from a neighborhood of it.
\end{longlist}
Let us start by evaluating the Laplace transform
of $V_{\omega}(x)+V_{\omega}(y)$.
We do it as a special case of a more general asymptotic formula
for the Laplace functional which we shall make repeated use of in
the sequel. We write $m_{\mu}$ for the
barycenter $\int x \mu(\d x)$ for a probability measure $\mu$
on $\R^d$.
%
%
\begin{proposition}
\label{Laplace}
Let $\varepsilon\in(0,1/\alpha)$.
For any probability measure $\mu$ supported inside
$B(0, t^{1/\alpha-\varepsilon})$,
\begin{eqnarray*}
&&
\E\bigl[\exp\bigl\{-t\langle\mu,V_{\omega} \rangle\bigr\}\bigr]
\\
&&\qquad =\exp\biggl\{-a_1t^{{d}/{\alpha}}- \bigl(C(d,\alpha)+o(1)
\bigr)t^{({d-2})/{\alpha}}\int|z-m_{\mu}|^2\mu(\d z) \biggr\}
\end{eqnarray*}
as $t\to\infty$.
In particular,
%
%
\begin{equation}
\label{2pt-Laplace} \qquad
\E\biggl[\exp\biggl\{-\frac{t}{2}
\bigl(V_{\omega}(x)+V_{\omega
}(y)\bigr) \biggr\} \biggr] \le\exp\bigl
\{-a_1t^{{d}/{\alpha}} -c_1t^{({d-2})/{\alpha}}|x-y|^2
\bigr\}
\end{equation}
for $c_1=C(d,\alpha)/5$ and $x,y \in\R^d$ with
$|x-y|<t^{1/\alpha-\varepsilon}$ when $t$ is sufficiently large.
\end{proposition}
%
%
\begin{pf}
We may assume $m_{\mu}=0$ by translation invariance of Poisson
point process. By a well-known formula for Laplace functional of
Poisson point process,
\begin{eqnarray*}
&&-{\log\E}\bigl[\exp\bigl\{-t\langle\mu,V_{\omega} \rangle\bigr\}\bigr]
\\
&&\qquad =\int\bigl(1-e^{-t\int\hat v(z-y) \mu(\d z)} \bigr)\,\d y
\\
&&\qquad =\int\bigl(1-e^{-t\hat v(-y)} \bigr)\,\d y +\int\bigl(e^{-t\hat v(-y)}-
e^{-t\int\hat v(z-y) \mu(\d z)} \bigr)\,\d y.
\end{eqnarray*}
This first term is easily shown to be $a_1t^{d/\alpha}+O(e^{-t})$;
see Lemma 1 of \cite{Fuk11a}.
We pick a positive constant
\[
0<\delta<\frac{\varepsilon}{\alpha+2}
\]
and divide the second term into the integrals over
$\{|y| <t^{1/\alpha-\delta}\}$ and
$\{|y|\ge t^{1/\alpha-\delta}\}$.
For the first region, we know
\[
\min\biggl\{t\hat v(-y), t\int\hat v(z-y) \mu(\d z) \biggr\} \ge
t^{\alpha\varepsilon}
\]
by the assumption on $\mu$, and hence the contribution form
this part is negligible.
For the second region, we may replace $\hat v$ by
$v(\cdot)=|\cdot|^{-\alpha}$ and
by a change of variable, it follows that
\begin{eqnarray*}
&&\int_{|y|\ge t^{1/\alpha-\delta}} \bigl(e^{-t\hat v(-y)}- e^{-t\int
\hat v(z-y) \mu(\d x)} \bigr)
\,\d y
\\
&&\qquad = t^{{d}/{\alpha}}\int_{|\eta|\ge t^{-\delta}} e^{-|\eta|^{-\alpha
}} \bigl(1-
e^{-\int|\eta-t^{-{1}/{\alpha}}z|^{-\alpha}-|\eta|^{-\alpha}
\mu(\d z)} \bigr)\,\d\eta.
\end{eqnarray*}
By Taylor's theorem, we can find bounded functions
$R_2$ and $R_3$ such that
\begin{eqnarray*}
&&\bigl|\eta-t^{-{1}/{\alpha}}z\bigr|^{-\alpha}-|\eta|^{-\alpha}
\\
&&\qquad =t^{-{1}/{\alpha}}\bigl\langle z, \nabla v(\eta)\bigr\rangle
+R_2(z,\eta)t^{-{2}/{\alpha}}|z|^2|\eta|^{-\alpha-2}
\\
&&\qquad =t^{-{1}/{\alpha}}\bigl\langle z, \nabla v(\eta)\bigr\rangle+
\tfrac{1}{2}t^{-{2}/{\alpha}}\bigl\langle z, \mathrm{Hess}_v(\eta
)z\bigr\rangle+R_3(z,\eta)t^{-{3}/{\alpha}}|z|^3|
\eta|^{-\alpha-3}
\end{eqnarray*}
for $|\eta|\ge t^{-\delta}$ and $|z| \le t^{1/\alpha-\varepsilon}$.
Using this second line and recalling that $m_{\mu}=0$,
we get
\[
\biggl\llvert\int\bigl|\eta-t^{-{1}/{\alpha}}z\bigr|^{-\alpha}-|\eta
|^{-\alpha}
\mu(\d z)\biggr\rrvert\le\|R_2\|_{\infty}O
\bigl(t^{-2\varepsilon+\delta(\alpha+2)}\bigr),
\]
when $|\eta|\ge t^{-\delta}$.
This right-hand side is $o(1)$, thanks to our choice of $\delta$,
and thus we can use the inequality $|1-e^{-a}-a|<a^2$ which holds for
small $a$ to obtain
%
%
\begin{eqnarray}
\label{expanded}
&&
\int_{|\eta|\ge t^{-\delta}} e^{-|\eta|^{-\alpha}}
\bigl(1-
e^{-\int|\eta-t^{-{1}/{\alpha}}z|^{-\alpha}-|\eta|^{-\alpha}
\mu(\d z)} \bigr)
\nonumber
\\
&&\qquad =\frac{1}{2}\int_{|\eta|\ge t^{-\delta}} e^{-|\eta|^{-\alpha
}}t^{-{2}/{\alpha}}
\int\bigl\langle z, \mathrm{Hess}_v(\eta)z\bigr\rangle\mu(\d z)\,\d\eta
\nonumber\\[-8pt]\\[-8pt]
&&\qquad\quad{} +\int_{|\eta|\ge t^{-\delta}} e^{-|\eta|^{-\alpha}} t^{-{3}/{\alpha}}\int
R_3(z,\eta)|z|^3\mu(\d z) |\eta|^{-\alpha-3}\,\d\eta
\nonumber
\\
&&\qquad\quad{} + O \biggl(\int_{|\eta|\ge t^{-\delta}} e^{-|\eta|^{-\alpha}}
\biggl(t^{-{2}/{\alpha}}\int R_2(z,\eta)|z|^2\mu(\d z)
\biggr)^2 |\eta|^{-2\alpha-4}\,\d\eta\biggr).\nonumber
\end{eqnarray}
A computation shows that the first term in the right-hand side
of (\ref{expanded}) equals
\begin{eqnarray*}
&&\frac{1}{2}t^{-{2}/{\alpha}} \int\biggl\langle z, \int
_{|\eta|\ge t^{-\delta}} \mathrm{Hess}_v(\eta)e^{-|\eta|^{-\alpha}}\,\d
\eta
z \biggr\rangle\mu(\d z)
\\
&&\qquad =\bigl(C(d,\alpha)+o(1)\bigr)t^{-{2}/{\alpha}} \int|z|^2\mu(\d z).
\end{eqnarray*}
The other terms in (\ref{expanded}) turn out to be smaller order
than this: indeed by the assumption on $\mu$:
\begin{longlist}
\item the second term is bounded by
$t^{-2/\alpha}\int|z|^2\mu(\d z)$ multiplied by
\[
O \bigl(t^{-{1}/{\alpha}} \bigr) \sup\bigl\{|z|\dvtx z \in
\operatorname{supp}\mu\bigr\}=o(1);
\]
\item$t^{-2/\alpha}\int|z|^2\mu(\d z)$
itself is of $o(1)$, and the third term has smaller order than it
as its square.\qed
\end{longlist}
\noqed\end{pf}
From the above lemma, we can deduce controls on
$V_{\omega}(x) + V_{\omega}(y)$ for all
$x,y$ within an intermediate distance.
%
%
\begin{lemma}
\label{lem4}
For any $\varepsilon\in(0,1/\alpha)$, there exists $R>0$ such that
%
%
\begin{eqnarray}
\label{two-point}
\qquad&&
\lim_{t \to\infty}Q_t \bigl(V_{\omega}(x)+V_{\omega}(y)
> 2h_t+c_1t^{-({\alpha-d+2})/{\alpha}}|x-y|^2
\nonumber\\
&&\hspace*{39.7pt}\mbox{ for all } x,y \in(-2t,2t)^d \mbox{ with
}Rt^{({\alpha-d+6})/({8\alpha})}<|x-y| <t^{1/\alpha-\varepsilon} \bigr
)\\
&&\qquad= 1.
\nonumber
\end{eqnarray}
\end{lemma}
%
%
\begin{pf}
We prove that for sufficiently large $R>0$,
%
%
\begin{eqnarray}
\label{two-point-decay}
&&
Q_t \bigl( V_{\omega}(x)+V_{\omega}(y)
\le2h_t+c_1t^{-({\alpha-d+2})/{\alpha}}|x-y|^2 \bigr)\nonumber\\[-8pt]\\[-8pt]
&&\qquad =
\O\bigl(\exp\bigl\{-t^{({\alpha+d-2})/({2\alpha})} \bigr\}
\bigr)\nonumber
\end{eqnarray}
uniformly in $x,y \in\R^d$ with
$Rt^{({\alpha-d+6})/({8\alpha})}<|x-y|<t^{1/\alpha-\varepsilon}$.
One can deduce (\ref{two-point}) from this by dividing
$(-2t,2t)^d$ into small cubes and using the union bound, just
as in the proof of Lemma~\ref{lem1}.

To prove (\ref{two-point-decay}), we may restrict our
consideration to the event
\[
E= \Bigl\{\min_{0\le s \le t} V_{\omega}(X_s)-h_t
\ge-2a_2t^{-({3\alpha-3d+2})/({4\alpha})} \Bigr\}
\]
by Lemma~\ref{lem2}.
Then denoting the event in (\ref{two-point-decay}) by $F_{x,y}$,
we have
\begin{eqnarray*}
&&\E\otimes E_0 \biggl[\exp\biggl\{-\int_0^{t}
V_{\omega}(X_s)\,\d s \biggr\}\dvtx E \cap F_{x,y} \biggr]
\\
&&\qquad \le\exp\biggl\{-a_1\frac{d}{\alpha}t^{{d}/{\alpha}}
+2a_2t^{({\alpha+3d-2})/({4\alpha})} \biggr\}\P(F_{x,y}).
\end{eqnarray*}
Now a simple large deviation bound with the help
of (\ref{2pt-Laplace}) shows
\[
\P(F_{x,y}) \le\exp\biggl\{-a_1\frac{\alpha-d}{\alpha}t^{{d}/{\alpha}}
-\frac{c_1}{2}R^2t^{({\alpha+3d-2})/({4\alpha})} \biggr\}
\]
uniformly in $x,y \in\R^d$ with
$Rt^{({\alpha-d+6})/({8\alpha})}<|x-y|<t^{1/\alpha-\varepsilon}$.
Therefore by taking sufficiently large $R$, we obtain
\[
Q_t(E \cap F_{x,y}) =\O\bigl(\exp\bigl
\{-t^{({\alpha+3d-2})/({4\alpha})} \bigr\} \bigr).
\]
Since $\frac{\alpha+3d-2}{4\alpha}>\frac{\alpha+d-2}{2\alpha}$
for $\alpha<d+2$, this completes the proof.
\end{pf}
Let $A_R(x)$ denote the annulus
$B (x,2Rt^{({\alpha-d+6})/({8\alpha})} )
\setminus B (x,Rt^{({\alpha-d+6})/({8\alpha})} )$
around $x \in\R^d$.
The above lemma implies that
%
%
\begin{eqnarray}
\label{annulus}
&&
\lim_{t \to\infty}Q_t \bigl(\mbox{for any }x \in
\L\cap(-t,t)^d,
\nonumber\\[-8pt]\\[-8pt]
&&\hspace*{40pt} V_{\omega}>h_t+c_1R^2t^{-({3\alpha-3d+2})/({4\alpha})}
\mbox{ on }A_R(x) \bigr)=1.\nonumber
\end{eqnarray}
As a consequence, if $c_1R^2>1$, then
every connected component of $\L\cap(-t,t)^d$ is contained in
a ball with radius
$Rt^{({\alpha-d+6})/({8\alpha})}$
with high probability.
Therefore Proposition~\ref{weak-localization} reduces to the
following.
%
%
\begin{proposition}
There exists $R>0$ such that
\[
\lim_{t \to\infty}Q_t \Bigl(\sup_{0 \le s \le t}
\operatorname{dist}(X_s, \L) \le2R t^{({\alpha-d+6})/({8\alpha})} \Bigr)=1.
\]
\end{proposition}
%
%
\begin{pf}
We use two lemmas whose proofs are given later.
The first one is a lower bound for $Z_t$ in terms of
a random eigenvalue, which is not explicit, but much more
precise than (\ref{Fuk10}).
%
%
\begin{lemma}
\label{lem5}
There exists $c_2>0$ such that
\[
Z_t \ge c_2t^{-2d}\E\bigl[\exp\bigl\{-t \lambda_1^{\omega}\bigl((-t,t)^d\bigr)\bigr\}
\bigr].
\]
\end{lemma}
%
The second one is an upper bound for the eigenvalue.
%
%
\begin{lemma}
\label{lem6}
For any $\varepsilon>0$,
\[
\lim_{t \to\infty}Q_t \bigl(\lambda_1^{\omega}
\bigl((-t,t)^d\bigr)\le h_t+ \varepsilon t^{-({3\alpha-3d+2})/({4\alpha})}
\bigr)= 1.
\]
\end{lemma}
%
Let $G$ be the intersection of $\{X_{[0,t]}\subset(-t,t)^d\}$,
$\{X_{[0,t]}\cap\L\neq\varnothing\}$ and the events in (\ref{annulus})
and in Lemma~\ref{lem6}.
As we have proved $\lim_{t\to\infty}Q_t(G)=1$ above, we restrict
ourselves to $G$. Though we drop $\cap G$ from the notation for simplicity,
it is assumed throughout the proof.
Now, let us introduce stopping times
\begin{eqnarray*}
H_0&=&\inf \bigl\{s \in[0,t]\dvtx\operatorname{dist}(X_s,
\L) > 2R t^{({\alpha-d+6})/({8\alpha})} \bigr\},
\\
H_1&=&\inf \bigl\{s \in[0,t]\dvtx X_s \in\L \bigr\}.
\end{eqnarray*}
We first consider the case $H_0 > H_1$.
In this case, we further introduce random times defined by
\begin{eqnarray*}
H_2&=&\inf \bigl\{s \in[H_1,t]\dvtx X_s
\notin B \bigl(X_{H_1},2Rt^{({\alpha-d+6})/({8\alpha})} \bigr) \bigr\},
\\
H_3&=&\sup \bigl\{s \in[H_1,H_2]\dvtx
X_s \in B \bigl(X_{H_1},Rt^{({\alpha-d+6})/({8\alpha})} \bigr) \bigr\},
\end{eqnarray*}
which satisfy $H_1<H_3<H_2<H_0<t$.
Let us define $E_k=\{H_3 \in[k,k+1)\}$,
$S=\{H_2-H_3>t^{({\alpha+d-2})/({2\alpha})}, H_2<t\}$ (slow crossing)
and $F=\{H_2-H_3 \le t^{({\alpha+d-2})/({2\alpha})}, H_2<t\}$
(fast crossing). Note that for $s \in[H_3,H_2]$,
we have $X_s \in A_R(X_{H_1})$ and hence
\begin{eqnarray*}
V_{\omega}(X_s) &\ge& h_t+c_1R^2
t^{-({3\alpha-3d+2})/({4\alpha})}
\\
&\ge&\lambda_1^{\omega}\bigl((-t,t)^d
\bigr)+t^{-({3\alpha-3d+2})/({4\alpha})}
\end{eqnarray*}
for sufficiently large $R>0$.
We use the strong Markov property at $H_2$ to see
%
%
\begin{eqnarray}
\label{crossing}
&& \E\otimes E_0 \biggl[\exp \biggl\{-\int
_0^{t} V_{\omega}(X_s)\,\d s
\biggr\}\dvtx S \cap E_k \biggr]
\nonumber
\\
&&\qquad \le\E \biggl[E_0 \biggl[\exp \biggl\{-\int
_0^{k+1} V_{\omega
}(X_s)\,\d
s \biggr\}
\nonumber\\
&&\hspace*{28.6pt}\qquad\quad{}\times \exp \bigl\{-(H_2-k-1)\\
&&\hspace*{100.5pt}{}\times
\bigl(\lambda_1^{\omega}
\bigl((-t,t)^d\bigr)+t^{-({3\alpha-3d+2})/({4\alpha})} \bigr) \bigr\}
1_{S\cap E_k} \nonumber
\\
&&\hspace*{28.6pt}\qquad\quad{}\times  E_{X_{H_2}} \biggl[ \exp \biggl\{-\int_0^{t-H_2}V_{\omega}(X_s)
\,\d s \biggr\}\dvtx X_{[0, t-H_2]}\subset(-t,t)^d \biggr] \biggr]
\biggr].\nonumber
\end{eqnarray}
Now from a uniform bound
%
%
\begin{eqnarray}
\label{unif-FK} && \sup_{x \in\R^d}E_x \biggl[\exp \biggl\{-
\int_0^{s} V_{\omega
}(X_s)
\,\d s \biggr\}\dvtx X_{[0,s]}\subset(-t,t)^d \biggr]
\nonumber\\[-8pt]\\[-8pt]
&&\qquad \le c(d) \bigl(1+\bigl(s\lambda_1^{\omega}
\bigl((-t,t)^d\bigr)\bigr)^{d/2} \bigr)\exp\bigl\{-s
\lambda_1^{\omega}\bigl((-t,t)^d\bigr)\bigr\}\nonumber
\end{eqnarray}
on Feynman--Kac semigroup (see Theorem 3.1.2
in \cite{Szn98})
and the fact $H_2-k-1>t^{({\alpha+d-2})/({2\alpha})}-1$ on $S\cap E_k$,
it follows that
%
%
\begin{equation}
\label{slow}\qquad  \frac{\E\otimes E_0 [\exp\{-\int_0^{t} V_{\omega
}(X_s)\,\d s \}\dvtx S \cap E_k ]} {
\E[\exp\{-t \lambda_1^{\omega}((-t,t)^d)\} ]} =\O \bigl(\exp \bigl\{-t^{({d+2-\alpha})/({4\alpha})} \bigr\}
\bigr).
\end{equation}
As for the fast crossing, we proceed as in (\ref{crossing}) and
use the Markov property at time $k$ to get
\begin{eqnarray*}
\hspace*{-4pt}&& \E\otimes E_0 \biggl[\exp \biggl\{-\int_0^{t}
V_{\omega}(X_s)\,\d s \biggr\}\dvtx F \cap E_k
\biggr]
\\
\hspace*{-4pt}&&\hspace*{-2pt}\qquad \le\E \biggl[E_0 \biggl[\exp \biggl\{-\int
_0^{k} V_{\omega
}(X_s)\,\d
s \biggr\}\exp \bigl\{-(H_2-k-1) \lambda_1^{\omega}
\bigl((-t,t)^d\bigr) \bigr\}1_{E_k}
\\
\hspace*{-4pt}&&\hspace*{-2pt}\qquad\quad\hspace*{28pt}{}\times P_{X_k} \bigl( \sup \bigl\{|X_s-X_0|
\dvtx0 \le s \le t^{({\alpha+d-2})/({2\alpha})} \bigr\} >Rt^{({\alpha-d+6})/({8\alpha})} \bigr)
\\
\hspace*{-4pt}&&\hspace*{-2pt}\qquad\quad\hspace*{55pt}{}\times E_{X_{H_2}} \biggl[ \exp \biggl\{-\int_0^{t-H_2}V_{\omega}(X_s)
\,\d s \biggr\}\dvtx X_{[0, t-H_2]}\subset(-t,t)^d \biggr] \biggr]
\biggr].
\end{eqnarray*}
Since we have
%
%
\begin{eqnarray}
\label{fast} && P_{X_k} \bigl( \sup \bigl\{|X_s-X_0|
\dvtx 0\le s \le t^{({\alpha+d-2})/({2\alpha})} \bigr\} >Rt^{({\alpha-d+6})/({8\alpha})} \bigr)
\nonumber
\\
&&\qquad =P_0 \bigl( \sup \bigl\{|X_s|\dvtx0\le s
\le t^{({\alpha+d-2})/({2\alpha})} \bigr\} >Rt^{({\alpha-d+6})/({8\alpha})} \bigr)
\\
&&\qquad =\O \bigl(\exp \bigl\{-t^{({d+2-\alpha})/({4\alpha})} \bigr\}
\bigr),\nonumber
\end{eqnarray}
it follows that
\[
\frac{\E\otimes E_0 [\exp\{-\int_0^{t} V_{\omega
}(X_s)\,\d s \}\dvtx F \cap E_k ]} {
\E[\exp\{-t \lambda_1^{\omega}((-t,t)^d)\} ]} =\O \bigl(\exp \bigl\{-t^{({d+2-\alpha})/({4\alpha})} \bigr\} \bigr)
\]
just as before.
Summing (\ref{slow}) and (\ref{fast}) over $k \in[0,t]\cap\Z$
and using Lemma~\ref{lem5}, we obtain
%
%
\begin{equation}
\lim_{t\to\infty} Q_t \Bigl(\sup_{s\in[H_1, t]}
|X_s|>2Rt^{({\alpha-d+6})/({8\alpha})} \Bigr)=0.
\end{equation}

In the other case $H_0< H_1$, we introduce random times defined by
\begin{eqnarray*}
\tilde{H}_2&=&\sup \Bigl\{s \in[0,t]\dvtx\sup_{0 \le s \le t}
\operatorname{dist}(X_s, \L)>2Rt^{({\alpha-d+6})/({8\alpha})} \Bigr\},
\\
\tilde{H}_3&=&\inf \bigl\{s \in[\tilde{H}_2,t]\dvtx
\operatorname{dist}(X_s, \L) \le Rt^{({\alpha-d+6})/({8\alpha})} \bigr\},
\end{eqnarray*}
which satisfy $H_0<\tilde{H}_2<\tilde{H}_3<H_1<t$.
By using $\tilde{H}_2$ and $\tilde{H}_3$
instead of $H_2$ and~$H_3$, respectively, one can show that
\[
\lim_{t\to\infty} Q_t \Bigl(\sup_{s\in[0, H_1]}
|X_s|>2Rt^{({\alpha-d+6})/({8\alpha})} \Bigr)=0
\]
by the same argument as above.
\end{pf}
%
%
\begin{pf*}{Proof of Lemma~\ref{lem5}}
We start by introducing the notation
\[
T_t^{\omega}f(x) =E_x \biggl[f(X_t)
\exp \biggl\{-\int_0^{t} V_{\omega}(X_s)
\,\d s \biggr\}\dvtx X_{[0,t]}\subset(-t,t)^d \biggr].
\]
Since we know from (\ref{macro-box})
that
\[
Z_t=\E \bigl[T_t^{\omega}1(0)
\bigr]+o(Z_t),
\]
we consider the first term on the right-hand side.

Using translation invariance with respect to $\P$, we find
%
%
\begin{eqnarray}
\label{102} Z_t &=&\E\otimes E_0 \biggl[ \exp \biggl
\{- \int_0^{t} V_{\omega}(X_s)
\,\d s \biggr\} \biggr]
\nonumber
\\
&=&\frac{1}{(2t)^d}\E \biggl[ \int_{(-t,t)^d} E_x
\biggl[ \exp \biggl\{-\int_0^{t}
V_{\omega
}(X_s)\,\d s \biggr\} \biggr]\,\d x \biggr]
\nonumber
\\
&\ge&\frac{1}{(2t)^d}\E \biggl[ \int_{(-t,t)^d}
T_t^{\omega}1(x) \,\d x\dvtx\lambda_1^{\omega}
\bigl((-t,t)^d\bigr)\le1 \biggr]
\\
&=& \frac{1}{(2t)^d}\E \bigl[ \bigl\langle T_t^{\omega}1,1
\bigr\rangle\dvtx\lambda_1^{\omega}\bigl((-t,t)^d
\bigr)\le1 \bigr]
\nonumber
\\
&\ge&\frac{1}{(2t)^d}\E \bigl[ \bigl\langle\phi^{\omega}_1,1
\bigr\rangle^2 \exp\bigl\{-t\lambda_1^{\omega}
\bigl((-t,t)^d\bigr)\bigr\}\dvtx\lambda_1^{\omega}
\bigl((-t,t)^d\bigr)\le1 \bigr],\nonumber
\end{eqnarray}
where $\phi_1^{\omega}$ is the $L^2$-normalized nonnegative
eigenfunction associated with $\lambda_1^{\omega}((-t, t)^d)$.
Since $\|T_1^{\omega}\|_{1 \to\infty}\le(2\pi)^{-{d}/{2}}$, we have
\[
0 \le\exp\bigl\{-\lambda_1^{\omega}\bigl((-t,t)^d
\bigr)\bigr\} \phi_1^{\omega}(z) = T_1^{\omega}
\phi_1^{\omega}(z) \le(2\pi)^{-{d}/{2}} \bigl\langle
\phi_1^{\omega},1 \bigr\rangle
\]
for any $z \in(-t,t)^d$. Therefore on $\{\lambda_1^{\omega}((-t,t)^d)\le1\}$,
\[
\bigl\langle\phi_1^{\omega},1 \bigr\rangle^2 \ge(2
\pi)^d e^{-2} \sup_{z\in(-t,t)^d}\phi_1^{\omega}(z)^2
\ge(2\pi)^d e^{-2} (2t)^{-d},
\]
where we have used $\operatorname{supp} \phi_1^{\omega}\subset(-t,t)^d$
together with $\|\phi_1^{\omega}\|_2=1$
in the second inequality. Coming back to (\ref{102}), we obtain
\[
Z_t \ge\frac{(2\pi)^d e^{-2}}{(2t)^{2d}} \E \bigl[\exp \bigl\{-t
\lambda_1^{\omega}\bigl((-t,t)^d\bigr)\bigr\}\dvtx
\lambda_1^{\omega}\bigl((-t,t)^d\bigr)\le1 \bigr]
\]
for $t \ge1$.
We can drop $\{\lambda_1^{\omega}((-t,t)^d)\le1\}$ from the right-hand side since
\[
\E \bigl[\exp\bigl\{-t\lambda_1^{\omega}\bigl((-t,t)^d
\bigr)\bigr\}\dvtx\lambda_1^{\omega}\bigl((-t,t)^d
\bigr)> 1 \bigr]\le e^{-t}=o(Z_t).
\]
\upqed\end{pf*}
%
%
\begin{pf*}{Proof of Lemma~\ref{lem6}}
The argument is similar to that for Lemma~\ref{lem2}.
We may restrict ourselves to $\{X_{[0,t]}\subset(-t, t)^d\}$
thanks to (\ref{macro-box}).
Then by (\ref{unif-FK}), it follows that
\begin{eqnarray*}
&&\E\otimes E_0 \biggl[\exp \biggl\{-\int_0^{t}
V_{\omega}(X_s)\,\d s \biggr\}\dvtx X_{[0,t]}
\subset(-t, t)^d, \lambda_1^{\omega}
\bigl((-t,t)^d\bigr)> 3a_1t^{-({\alpha-d})/{\alpha}} \biggr]
\\
&&\qquad \le\exp \bigl\{-2a_1t^{{d}/{\alpha}} \bigr\} c(d)
\sup_{\lambda>0} \bigl(1+\lambda^{d/2} \bigr)e^{-{\lambda}/{3}}
\\
&&\qquad =o(Z_t).
\end{eqnarray*}
Next we show that
%
%
\begin{eqnarray}
\label{unif-exp} &&\E\otimes E_0 \biggl[\exp \biggl\{-\int
_0^{t} V_{\omega}(X_s)\,\d s
\biggr\}\dvtx X_{[0,t]}\subset(-t, t)^d,
\nonumber
\\
&&\hspace*{39.1pt}\lambda_1^{\omega}\bigl((-t,t)^d
\bigr)-h_t-\varepsilon t^{-({3\alpha-3d+2})/({4\alpha})} \in t^{-1}[k-1,k) \biggr]
\\
&&\qquad =Z_t\O \bigl(\exp \bigl\{-t^{({\alpha+d-2})/({2\alpha})} \bigr\}
\bigr)\nonumber
\end{eqnarray}
uniformly in $k \in\N\cap[0, 3a_1 t^{d/\alpha}+1]$.
Using (\ref{unif-FK}) again, we can bound the above left-hand side by
%
%
\begin{eqnarray}
\label{splitted} &&\E\bigl[c(d) \bigl(1+\bigl(t\lambda_1^{\omega}
\bigl((-t,t)^d\bigr)\bigr)^{d/2} \bigr)\exp\bigl\{-t
\lambda_1^{\omega}\bigl((-t,t)^d\bigr)\bigr\} \dvtx
\nonumber
\\
&&\hspace*{12pt}\lambda_1^{\omega}\bigl((-t,t)^d
\bigr)-h_t-\varepsilon t^{-({3\alpha-3d+2})/({4\alpha})} \in t^{-1}[k-1,k) \bigr]
\nonumber\\[-8pt]\\[-8pt]
&&\qquad \le t^d \exp \bigl\{-th_t- \varepsilon
t^{({\alpha+3d-2})/({4\alpha})}-k+1 \bigr\}
\nonumber
\\
&&\qquad\quad\hspace*{0pt}{}\times\P \bigl(\lambda_1^{\omega}\bigl((-t,t)^d\bigr)
\le h_t+\varepsilon t^{-({3\alpha-3d+2})/({4\alpha})}+t^{-1}k
\bigr)\nonumber
\end{eqnarray}
for sufficiently large $t$.
To estimate the last probability, we use an inequality
\[
\P\bigl(\lambda_1^{\omega}\bigl((-t,t)^d\bigr)\le
\lambda\bigr) \le(2t)^d N(\lambda),
\]
where $N(\lambda)$ is the integrated density of states of
$-1/2\Delta+V_{\omega}$; see, for example,~\cite{CL90}, Chapter VI, for
the definition of integrated density of states and also the above
inequality. The asymptotics of $N(\lambda)$ has been determined up to
the second term in Theorem 3 of \cite{Fuk11a}:
\[
N(\lambda)=\exp \bigl\{-l_1 \lambda^{-{d}/({\alpha-d})} -
\bigl(l_2+o(1) \bigr) \lambda^{-({\alpha+d-2})/({2(\alpha-d)})} \bigr\}
\]
as $\lambda\downarrow0$, where
\[
l_1=\frac{\alpha-d}{\alpha} \biggl(\frac{d}{\alpha} \biggr)^{{d}/({\alpha-d})}
a_1^{{\alpha}/({\alpha-d})},\qquad l_2=a_2 \biggl(
\frac{da_1}{\alpha} \biggr)^{({\alpha+d-2})/({2(\alpha-d)})}.
\]
Substituting $h_t+\varepsilon t^{-({3\alpha-3d+2})/({4\alpha})}+t^{-1}k$
into $\lambda$, one can find a constant $c_3>0$ satisfying
\begin{eqnarray*}
&&\P \bigl(\lambda_1^{\omega}\bigl((-t,t)^d\bigr)
\le h_t+\varepsilon t^{-({3\alpha-3d+2})/({4\alpha})}+t^{-1}k \bigr)
\\
&&\qquad \le\exp \biggl\{-a_1\frac{\alpha-d}{\alpha}t^{{d}/{\alpha}} +
\varepsilon t^{({\alpha+3d-2})/({4\alpha})}+k - \bigl(a_2+c_3
\varepsilon^2 \bigr)t^{({\alpha+d-2})/({2\alpha})} \biggr\}
\end{eqnarray*}
for all $k \in\N\cap[0, 3a_1 t^{d/\alpha}+1]$
by a straightforward calculation. (It suffices to
take $c_3$ so small that
\[
(1+x)^{-{d}/({\alpha-d})}\ge1-\frac{d}{\alpha-d}x+c_3x^2
\]
for all $x \in[0, 4\alpha/d]$.)
From this and (\ref{splitted}), the desired
bound (\ref{unif-exp}) follows.
Summing over $k$, we obtain
\begin{eqnarray*}
&&\E\otimes E_0 \biggl[\exp \biggl\{-\int_0^{t}
V_{\omega}(X_s)\,\d s \biggr\}\dvtx
X_{[0,t]}
\subset(-t, t)^d,
\\
&&\qquad\hspace*{17.2pt}\lambda_1^{\omega}\bigl((-t,t)^d\bigr)
\in\bigl(h_t+\varepsilon t^{-({3\alpha-3d+2})/({4\alpha})}, 3a_1t^{-({\alpha-d})/{d}}
\bigr] \biggr]
\\
&&\qquad =Z_t\O \bigl(\exp \bigl\{-t^{({\alpha+d-2})/({2\alpha})} \bigr\}
\bigr),
\end{eqnarray*}
and we are done.
\end{pf*}
%
As a corollary to the weak localization, we obtain an upper
bound on the variance of $L_t$. In contrast to
Proposition~\ref{weak-localization}, the bound is of correct
order, and this will be crucial in the next section.
%
%
\begin{corollary}
\label{variance-bound}
For any $M_2 \ge\frac{2a_2}{C(d,\alpha)}$,
\[
\lim_{t\to\infty} Q_t \biggl(\int|x-m_{L_t}|^2L_t(
\d x) <M_2 t^{({\alpha-d+2})/({2\alpha})} \biggr) =1.
\]
\end{corollary}
%
%
\begin{pf}
In view of Proposition~\ref{weak-localization}, we only need to show
\begin{eqnarray*}
&&\lim_{t\to\infty} Q_t \biggl(\int|x-m_{L_t}|^2
L_t(\d x) \ge\frac{2a_2}{C(d,\alpha)} t^{({\alpha-d+2})/({2\alpha})},
\\
&&\hspace*{98pt}\operatorname{supp}L_t \subset B \bigl(0,M_1t^{({\alpha-d+6})/({8\alpha})}
\bigr) \biggr) =0.
\end{eqnarray*}
But it follows from Proposition~\ref{Laplace}
that if $L_t$ lies in the above event,
\begin{eqnarray*}
&&\E \biggl[\exp \biggl\{-\int_0^{t}
V_{\omega}(X_s)\,\d s \biggr\} \biggr]
\\
&&\qquad =\E \bigl[\exp \bigl\{-t\langle L_t, V_{\omega}
\rangle \bigr\} \bigr]
\\
&&\qquad =\exp \biggl\{-a_1t^{{d}/{\alpha}}- \bigl(C(d,\alpha)+o(1)
\bigr)t^{({d-2})/{\alpha}}\int|x-m_{L_t}|^2L_t(\d
x) \biggr\}
\\
&&\qquad \le\exp \bigl\{-a_1t^{{d}/{\alpha}}- \bigl(2a_2+o(1)
\bigr)t^{({\alpha+d-2})/({2\alpha})} \bigr\}
\\
&&\qquad =o(Z_t).
\end{eqnarray*}
\upqed\end{pf}
%
\section{Potential confinement}
\label{weak-to-PC}
We prove the following version of
Theorem~\ref{potential-confinement} in this section
by using the weak localization result.
It will be used in next section to
complete the proof of Theorem~\ref{localization}.
The proof of Theorem~\ref{potential-confinement} will be
completed at the end of next section with the help of
Theorem~\ref{localization}.
%
%
\begin{proposition}
\label{PC-L_t}
There exists an $\varepsilon>0$ such that
\begin{eqnarray*}
&& \lim_{t\to\infty} Q_t
\bigl(\sup \bigl\{\bigl|V_{\omega}(x)-V_{\omega}(m_{L_t})
-p_t(x-m_{L_t})\bigr|\dvtx
\\
&&\hspace*{78.2pt}|x-m_{L_t}|\le2t^{({\alpha-d+2})/({4\alpha})} \log t \bigr\}<t^{-
({\alpha-d+2})/({2\alpha})-\varepsilon}
\bigr)\\
&&\qquad= 1.
\end{eqnarray*}
\end{proposition}
We define the weighted measure by
\[
\frac{\d\P_t^{\mu}}{\d\P} =\frac{e^{-t\langle\mu,V_{\omega}\rangle}} {
\E[e^{-t\langle\mu, V_{\omega}\rangle}]}
\]
for each probability measure $\mu$ on $\R^d$.
Under $\P_t^{\mu}$, $\omega$ is a Poisson point process
with intensity
$e^{-t\int\hat v(x-y)\mu(\d x)}\,\d y$.
We take $M_2$ as in Corollary~\ref{variance-bound} and
define a class of probability measures on $\R^d$ by
\begin{eqnarray*}
\mathcal{P}_{\mathrm{wl}}&=& \biggl\{\mu\dvtx\operatorname{supp} \mu \subset B
\bigl(0,M_1t^{({\alpha-d+6})/({8\alpha})} \bigr),\\
&&\hspace*{6.5pt} \int|x-m_{\mu}|^2
\mu(\d x) <M_2 t^{-({\alpha-d+2})/({2\alpha})} \biggr\},
\end{eqnarray*}
where ``$\mathrm{wl}$'' stands for ``weakly localized.''
%
%
\begin{proposition}
\label{PC}
There exists an $\varepsilon>0$ such that uniformly in
$\mu\in\mathcal{P}_{\mathrm{wl}}$,
\begin{eqnarray*}
&& \lim_{t \to\infty} \P_t^{\mu} \bigl(\sup
\bigl\{\bigl|V_{\omega}(x)-V_{\omega}(m_{\mu}) -p_t(x-m_{\mu})\bigr|
\dvtx
\\
&&\hspace*{74.4pt}|x-m_{\mu}|\le2t^{({\alpha-d+2})/({4\alpha})} \log t \bigr
\}<t^{-
({\alpha-d+2})/({2\alpha})-\varepsilon} \bigr)\\
&&\qquad= 1.
\end{eqnarray*}
\end{proposition}
%
Let us first see that this immediately implies
Proposition~\ref{PC-L_t}.
Indeed, denoting the event in Proposition~\ref{PC} by $\mathtt{PC}_t(\mu)$,
we obtain
\begin{eqnarray*}
Q_t \bigl(\mathtt{PC}_t(L_t) \bigr) &\ge&
Q_t \bigl(\mathtt{PC}_t(L_t),
L_t\in \mathcal{P}_{\mathrm{wl}} \bigr)
\\
&=&\frac{1}{Z_t} E_0 \bigl[\E \bigl[e^{-t\langle L_t,V_{\omega}\rangle
} \bigr]
\P_t^{L_t} \bigl(\mathtt{PC}_t(L_t)
\bigr)\dvtx L_t\in\mathcal{P}_{\mathrm{wl}} \bigr]
\\
&\sim&\frac{1}{Z_t} E_0 \bigl[\E \bigl[e^{-t\langle L_t,V_{\omega}\rangle
}
\bigr] \dvtx L_t \in\mathcal{P}_{\mathrm{wl}} \bigr]
\\
& =& Q_t (L_t\in\mathcal{P}_{\mathrm{wl}} )\to1
\end{eqnarray*}
as $t\to\infty$ by
Proposition~\ref{weak-localization} and
Corollary~\ref{variance-bound}.%

To prove Proposition~\ref{PC}, we first
compute the expectation of $V_{\omega}(x)-V_{\omega}(m_{\mu})$
and then bound the variance.
%
%
\begin{lemma}
\label{lem7}
There exists an $\varepsilon>0$ such that uniformly in
$x \in\break B(0, 2t^{({\alpha-d+2})/({4\alpha})}\log t)$ and
$\mu\in\mathcal{P}_{\mathrm{wl}}$,
\[
\E_t^{\mu} \bigl[V_{\omega}(x)-V_{\omega}(m_{\mu})
\bigr] =p_t(x-m_{\mu})+ o \bigl(t^{-({\alpha-d+2})/({2\alpha
})-\varepsilon} \bigr)
\]
as $t\to\infty$.
\end{lemma}
%
%
\begin{pf}
By a well-known formula for Poisson point process,
\begin{eqnarray*}
&&\E_t^{\mu} \bigl[V_{\omega}(x)-V_{\omega}(m_{\mu})
\bigr]
\\
&&\qquad =\int \bigl(\hat v(x-y)-\hat v(m_{\mu}-y) \bigr)
e^{-t\int\hat v(z-y)\mu
(\d z)}\,\d y.
\end{eqnarray*}
We assume $m_{\mu}=0$ by translation.
Pick $\delta>0$ so small that
%
%
\begin{equation}
\label{delta} \frac{\alpha-d+6}{8\alpha}<\frac{1}{\alpha}-\delta
\quad\mbox{and}\quad \delta(
\alpha+2)<\frac{d+2-\alpha}{2\alpha}.
\end{equation}
Then uniformly in $\mu\in\mathcal{P}_{\mathrm{wl}}$,
\begin{eqnarray*}
&& \biggl\llvert\int_{|y|\le t^{{1}/{\alpha}-\delta}} \bigl(\hat v(x-y)-\hat v(-y)
\bigr) e^{-t\int\hat v(z-y)\mu(\d z)}\,\d y \biggr\rrvert
\\
&&\qquad \le\exp \bigl\{-t\inf \bigl\{|z-y|^{-\alpha}\dvtx|z|\le
M_1t^{({\alpha-d+6})/({8\alpha})}, |y|\le t^{{1}/{\alpha}-\delta} \bigr\}+O(\log t) \bigr\}
\\
&&\qquad =\O \bigl(\exp \bigl\{-t^{\alpha\delta} \bigr\} \bigr),
\end{eqnarray*}
and hence this part is negligible.
On the other hand, if $|y|>t^{1/\alpha-\delta}$, we can replace
$\hat{v}$ in the integrand by $v(\cdot)=|\cdot|^{-\alpha}$ and hence
\begin{eqnarray*}
&&\int_{|y|> t^{{1}/{\alpha}-\delta}} \bigl(\hat v(x-y)-\hat v(-y)
\bigr)e^{-t\int\hat v(z-y)\mu(\d z)}\,\d y
\\
&&\qquad =t^{-({\alpha-d})/{\alpha}}\int_{|\eta|> t^{-\delta}} \bigl(\bigl|\eta
-t^{-{1}/{\alpha}}x\bigr|^{-\alpha}-| \eta|^{-\alpha} \bigr) e^{-\int|\eta-t^{-{1}/{\alpha}}z|^{-\alpha}\mu
(\d z)}
\,\d\eta.
\end{eqnarray*}
Let us write
$D_t(\eta,x)=|\eta-t^{-{1}/{\alpha}}x|^{-\alpha}-|\eta
|^{-\alpha}$
for short and further decompose the integral in the last line as
%
%
\begin{eqnarray}
\label{expectation-relevant}
&&
t^{-({\alpha-d})/{\alpha}}\int_{|\eta|>
t^{-\delta}}
D_t(\eta,x)e^{-|\eta|^{-\alpha}}\nonumber\\[-8pt]\\[-8pt]
&&\qquad\hspace*{56pt}{} + D_t(\eta,x)
\bigl(e^{-\int D_t(\eta,z)\mu(\d z)}-1 \bigr)e^{-|\eta|^{-\alpha}} \,\d
\eta.\nonumber
\end{eqnarray}
Let us first see that the first term
of (\ref{expectation-relevant}) gives us the desired quantity.
We use the following expansions of $D_t(\eta,\cdot)$,
which immediately follow from
Taylor's theorem.
%
%
\begin{lemma}
\label{Taylor}
For $|z| \le M_1t^{({\alpha-d+6})/({8\alpha})}$ and
$|\eta|>t^{-\delta}$, there exist bounded functions $R_1$,
$R_2$ and $R_3$ such that
%
%
\begin{eqnarray}
\quad&& D_t(\eta,z)
\nonumber
\\
\label{Taylor1} &&\qquad=R_1(z,\eta)t^{-{1}/{\alpha}}|z||
\eta|^{-\alpha-1}
\\
\label{Taylor2} &&\qquad=t^{-{1}/{\alpha}} \bigl\langle z, \nabla v(\eta) \bigr
\rangle+R_2(z, \eta)t^{-{2}/{\alpha}}|z|^2|
\eta|^{-\alpha-2}
\\
\label{Taylor3} &&\qquad=t^{-{1}/{\alpha}} \bigl\langle z, \nabla v(\eta) \bigr
\rangle+ \tfrac{1}{2}t^{-{2}/{\alpha}} \bigl\langle z, \mathrm{Hess}_v(
\eta )z \bigr\rangle+R_3(z,\eta)t^{-{3}/{\alpha}}|z|^3|
\eta|^{-\alpha-3}.
\end{eqnarray}
\end{lemma}
%
Using (\ref{Taylor3}), one can deduce
\begin{eqnarray*}
&&t^{-({\alpha-d})/{\alpha}} \int_{|\eta|> t^{-\delta}} D_t(\eta,x) \exp
\bigl\{-|\eta|^{-\alpha} \bigr\}\,\d\eta
\\
&&\qquad =t^{-({\alpha-d+2})/{\alpha}}C(d,\alpha)|x|^2 +O \bigl(t^{-
({\alpha-d+2})/({2\alpha})-\varepsilon}
\bigr)
\end{eqnarray*}
for some $\varepsilon>0$ by a straightforward calculation.
[Note that $\int\nabla v(\eta) e^{-|\eta|^{-\alpha}}\,\d\eta=0$
by symmetry.]

Now let us turn to the estimate of the second term
of (\ref{expectation-relevant}).
As a consequence of (\ref{Taylor2}), we have
%
%
\begin{eqnarray}
\label{intensity-replacement} \biggl\llvert\int D_t(\eta,z) \mu(\d
z) \biggr\rrvert&=& \|R_2\|_{\infty}|\eta|^{-\alpha-2}t^{-{2}/{\alpha}}
\int|z|^2 \mu(\d x)
\nonumber\\[-8pt]\\[-8pt]
&=& O \bigl(t^{\delta(\alpha+2)-({d+2-\alpha})/({2\alpha})}
\bigr)\nonumber
\end{eqnarray}
uniformly in
$|\eta|>t^{-\delta}$ and $\mu\in\mathcal{P}_{\mathrm{wl}}$,
where we have used $m_{\mu}=0$.
Thanks to our choice of $\delta$, the last line goes to 0 as
$t \to\infty$.
Therefore we may use an inequality
$|e^{-a}-1| \le2|a|$ valid for small $a$ in the
second term of (\ref{expectation-relevant}) to obtain
\begin{eqnarray*}
&&t^{-({\alpha-d})/{\alpha}} \biggl\llvert\int_{|\eta|> t^{-\delta}}
D_t( \eta,x) \bigl(e^{-\int D_t(\eta,z)\mu(\d z)}-1 \bigr)e^{-|\eta|^{-\alpha
}} \,\d
\eta \biggr\rrvert
\\
&&\qquad \le2t^{-({\alpha-d})/{\alpha}} \int_{|\eta|> t^{-\delta}} \bigl|D_t(
\eta,x)\bigr| \biggl\llvert\int D_t(\eta,z)\mu(\d z) \biggr\rrvert
e^{-|\eta|^{-\alpha}} \,\d\eta.
\end{eqnarray*}
Applying (\ref{Taylor1}) and (\ref{Taylor2}) to $D_t(\eta,x)$
and $D_t(\eta,z)$, respectively, and then using the variance bound
for $\mu$, one can easily conclude that the above right-hand side
is of order $O(t^{-({\alpha-d+2})/({2\alpha})-\varepsilon})$
for some $\varepsilon>0$.
\end{pf}
\begin{pf*}{Proof of Proposition~\ref{PC}}
Throughout the proof we assume $\mu\in\mathcal{P}_{\mathrm{wl}}$,
and all the estimates below are supposed to be uniform in $\mu$.
We further assume \mbox{$m_{\mu}=0$} for simplicity.
That this causes no loss of generality will become clear in the
argument below.
In view of Lemma~\ref{lem7}, it suffices to estimate
the maximum of
\begin{eqnarray*}
&& \bigl\llvert V_{\omega}(x)-V_{\omega}(0) -\E_t^{\mu}
\bigl[V_{\omega}(x)-V_{\omega}(0) \bigr] \bigr\rrvert
\\
&&\qquad =\int \bigl(\hat v(x-y)-\hat v(-y) \bigr) \bigl(\omega(\d
y)-e^{-t\int
\hat v(z-y)\mu(\d z)}\,\d y \bigr)
\end{eqnarray*}
among $|x|\le r(t) \log t$
[recall $r(t)=t^{({\alpha-d+2})/({4\alpha})}$].
Fix $\delta>0$ as in (\ref{delta}).
We first show that the region $\{|y| \le t^{1/\alpha-\delta}\}$
makes only a negligible contribution.
Let us introduce the notation
\[
\bar{\omega}_t^{\mu}(\d y) =\omega(\d y)-
e^{-\int\hat v(z-y)\mu(\d z)} \,\d y.
\]
Observe that
\begin{eqnarray*}
&& \biggl\llvert\int_{|y| \le t^{{1}/{\alpha}-\delta}} \bigl(\hat v(x-y)-\hat v(-y)
\bigr) \bar{\omega}_t^{\mu}(\d y) \biggr\rrvert
\\
&&\qquad \le\int_{|y| \le t^{{1}/{\alpha}-\delta}} \bigl|\hat v(x-y)-\hat v(-y)\bigr| \bigl(
\omega(\d y)+e^{-t\int\hat v(z-y)\mu(\d z)}\,\d y \bigr)
\\
&&\qquad \le\int_{|y| \le t^{{1}/{\alpha}-\delta}} \bar{\omega}_t^{\mu}(
\d y) + 2\int_{|y| \le t^{{1}/{\alpha}-\delta}} e^{-t\int\hat v(z-y)\mu
(\d z)}\,\d y
\end{eqnarray*}
since $0 < \hat v \le1$.
The first term has zero $\E_t^{\mu}$-mean, and its variance,
which equals half of the second term, is seen to be of
$\O(\exp\{-t^{\alpha\delta}\})$
by the same argument as in the proof of Lemma~\ref{lem7}.
Therefore, Chebyshev's inequality yields
\[
\lim_{t\to\infty}\P_t^{\mu} \biggl(
\sup_{x \in\R^d} \biggl|\int_{|y| \le t^{{1}/{\alpha}-\delta}} \bigl(\hat v(x-y)-\hat v(-y)
\bigr) \bar{ \omega}_t^{\mu}(\d y) \biggr| <t^{-({\alpha-d+2})/({2\alpha})-\varepsilon}
\biggr)= 1
\]
for any $\varepsilon>0$.\vadjust{\goodbreak}

On the remaining part $\{|y|>t^{1/\alpha-\delta}\}$, we may
replace $\hat v$ by $v(\cdot)=|\cdot|^{-\alpha}$.
Then it follows by Lemma~\ref{Taylor} that
%
%
\begin{eqnarray}
\label{123} && \biggl\llvert\int_{|y| > t^{{1}/{\alpha}-\delta}} \bigl(\hat v(x-y)-
\hat v(-y) \bigr) \bar{\omega}_t^{\mu}(\d y) \biggr\rrvert
\nonumber
\\[2pt]
&&\qquad \le\biggl| \biggl\langle x, \int_{|y| > t^{{1}/{\alpha}-\delta}} \nabla v(y) \bar{
\omega}_t^{\mu}(\d y) \biggr\rangle\biggr|
\nonumber
\\[2pt]
&&\qquad\quad{}+ \biggl| \biggl\langle x, \int_{|y| > t^{{1}/{\alpha}-\delta}}
\mathrm{Hess}_v(y)\bar{\omega}_t^{\mu}(
\d y) x \biggr\rangle\biggr|
\\[2pt]
&&\qquad\quad{}+|x|^3 \biggl|\int_{|y| > t^{{1}/{\alpha}-\delta}}
R_3 \bigl(x,t^{-{1}/{\alpha}}y \bigr)|y|^{-\alpha-3} \bar{
\omega}_t^{\mu}(\d y) \biggr|
\nonumber
\\[2pt]
&&\qquad =:I_1(x)+I_2(x)+I_3(x).\nonumber
\end{eqnarray}
Let us start with $I_3(x)$. As we did for
$\{|y|\le t^{1/\alpha-\delta}\}$, we can bound the integral as
%
%
\begin{eqnarray}
\label{I_3} && \sup_{|x| \le r(t)\log t} \biggl\llvert\int
_{y>t^{{1}/{\alpha}-\delta}} R_3 \bigl(x,t^{-{1}/{\alpha}}y
\bigr)|y|^{-\alpha-3} \bar{\omega}_t^{\mu}(\d y) \biggr
\rrvert
\nonumber
\\[2pt]
&&\qquad\le\int_{|y| > t^{{1}/{\alpha}-\delta}} \|R_3\|_{\infty
}|y|^{-\alpha-3}
\bar{\omega}_t^{\mu}(\d y)
\\[2pt]
&&\qquad\quad{} + 2\int_{|y| > t^{{1}/{\alpha}-\delta}} \|R_3\|_{\infty
}|y|^{-\alpha-3}
e^{-\int|z-y|^{-\alpha}\mu(\d z)}\,\d y.\nonumber
\end{eqnarray}
The second term on the right-hand side can be seen to be of
$O(t^{-({\alpha-d+3})/{\alpha}})$,
by using the change of variable $y=t^{1/\alpha}\eta$
and then (\ref{intensity-replacement})
to replace\vspace*{1pt} $\int|\eta-t^{-1/\alpha}z|^{-\alpha}\mu(\d z)$ by
$|\eta|^{-\alpha}$.
Since
$|x|^3t^{-({\alpha-d+3})/{\alpha}}
=O (t^{-({\alpha-d+2})/({2\alpha})-\varepsilon} )$
uniformly in $|x| \le r(t)\log t$
for some $\varepsilon>0$, the second term is negligible.
By the same way, we can show that
%
%
\begin{eqnarray}
\label{var-I_3} && \var \biggl(\int_{|y| > t^{{1}/{\alpha}-\delta}}
\|R_3\|_{\infty}|y|^{-\alpha-3} \bar{\omega}_t^{\mu}(
\d y) \biggr)
\nonumber\\[-7pt]\\[-7pt]
&&\qquad =\int_{|y| > t^{{1}/{\alpha}-\delta}} \|R_3\|_{\infty
}^2|y|^{-2\alpha-6}
e^{-\int|z-y|^{-\alpha}\mu(\d z)}\,\d y\nonumber
\end{eqnarray}
is of order $O(t^{-({2\alpha-d+6})/{\alpha}})$.
Therefore by using Chebyshev's inequality, we obtain
\begin{eqnarray*}
&& \P_t^{\mu} \biggl( \sup_{|x| \le r(t)\log t}|x|^3
\biggl\llvert\int_{|y| > t^{{1}/{\alpha}-\delta}} \|R_3
\|_{\infty}|y|^{-\alpha-3} \bar{\omega}_t^{\mu}(
\d y) \biggr\rrvert\ge t^{-({\alpha-d+2})/({2\alpha})-\varepsilon} \biggr)
\\[2pt]
&&\qquad \le\P_t^{\mu} \biggl( \biggl\llvert\int
_{|y| > t^{{1}/{\alpha}-\delta}} \|R_3\|_{\infty}|y|^{-\alpha-3}
\bar{\omega}_t^{\mu}(\d y) \biggr\rrvert\ge
t^{-{5(\alpha
-d+2)}/({4\alpha})-{\varepsilon}/{2}} \biggr)
\\
&&\qquad =O \bigl(t^{-({3d+2-\alpha})/({2\alpha})+\varepsilon} \bigr),
\end{eqnarray*}
which goes to 0 as $t \to\infty$ for sufficiently small $\varepsilon$.

As for $I_1$ and $I_2$, we use variance bounds
\begin{eqnarray*}
\var \biggl(\int_{|y| > t^{{1}/{\alpha}-\delta}} \nabla v(y) \bar{
\omega}_t^{\mu}(\d y) \biggr) &=&O \bigl(t^{-({2\alpha-d+2})/{\alpha}}
\bigr),
\\
\var \biggl(\int_{|y| > t^{{1}/{\alpha}-\delta}} \mathrm{Hess}_v(y)
\bar{ \omega}_t^{\mu}(\d y) \biggr) &=&O
\bigl(t^{-({2\alpha-d+4})/{\alpha}} \bigr),
\end{eqnarray*}
which follow by routine arguments.
Given these bounds, one can easily conclude that
\[
\lim_{t\to\infty}\P_t^{\mu} \Bigl(
\sup_{|x| \le r(t)\log t} \bigl(I_1(x)+I_2(x) \bigr)
<t^{-({\alpha-d+2})/({2\alpha})-\varepsilon} \Bigr)= 1
\]
for some $\varepsilon>0$.
\end{pf*}
%
%
\begin{remark}
\label{rem1}
Inspecting the above proof, one can easily improve the
bound as
\begin{eqnarray*}
&& \lim_{t\to\infty} Q_t
\bigl(\sup \bigl\{\bigl|V_{\omega}(x)-V_{\omega}(m_{L_t})
-\E_t^{\mu} \bigl[V_{\omega}(x)-V_{\omega}(m_{L_t})
\bigr]\bigr|\dvtx
\\
&&\hspace*{122.5pt}|x-m_{L_t}|\le2t^{({\alpha-d+2})/({4\alpha})} \log t \bigr\}<t^{-\beta}
\bigr)\\
&&\qquad= 1
\end{eqnarray*}
for any
$\beta<\min\{\frac{3\alpha-d+2}{4\alpha},\frac{\alpha
-d+6}{4\alpha}\}$.
This will be used in Section~\ref{fluctuation}.
\end{remark}
%
\section{Strong localization}
\label{strong-localization}
In this section, we accomplish the proof of Theorem~\ref{localization}.
Loosely speaking, we do this by simply repeating the argument for
Proposition~\ref{weak-localization} but in the correct scale.
Recall that we used an annulus $A_R$ which was much larger than the
correct scale $r(t)=t^{({\alpha-d+2})/({4\alpha})}$ in the proof of
Proposition~\ref{weak-localization}. It was because we had control
of $V_{\omega}$ only at points in the intermediate distance from $\L
$; see Lemma~\ref{lem4}. Now that Proposition~\ref{PC-L_t} is available,
we have the control of $V_{\omega}$
near $m_{L_t}$ and can work in the correct scale.
\begin{pf*}{Proof of Theorem~\ref{localization}}
Set $G=\{L_t \in\mathcal{P}_{\mathrm{wl}}\}\cap\mathtt{PC}(L_t)$.
As we have already shown $\lim_{t\to\infty}Q_t(G)=1$, we restrict
ourselves to $G$. Though we drop $\cap G$ from the notation for
simplicity, it is assumed throughout the proof.

We first fix $q \in\Z^d$ and assume $m_{L_t} \in q+[0,1)^d$.
Let us set $R=\sqrt{2M_2}$ and define
\[
H_1(q)=\inf \bigl\{s \ge0\dvtx|X_s-q| \le
Rt^{({\alpha-d+2})/({4\alpha})} \bigr\}.\vadjust{\goodbreak}
\]
Clearly $\{H_1(q)<t\} \supset\{L_t \in\mathcal{P}_{\mathrm{wl}}\}
\cap\{m_{L_t} \in q+[0,1)^d\}$.
Note also that on $\mathtt{PC}(L_t)$, we have
%
%
\begin{equation}
\label{ev-bound} \lambda_1^{\omega}\bigl((-t,t)^d
\bigr)\le V_{\omega}(m_{L_t})+ \bigl(a_2+o(1)
\bigr)t^{-({\alpha-d+2})/({2\alpha})}
\end{equation}
by an application of the Rayleigh--Ritz variational formula
together with a simple cut-off argument.

Now, fix $\delta\in(0,1/2)$, and define
\[
H_2(q)=\inf \bigl\{s \ge0\dvtx|X_s-q|>
t^{({\alpha-d+2})/({4\alpha})}( \log t)^{{1}/{2}+\delta} \bigr\}.
\]
We first consider the case $H_1(q)<H_2(q)\le t$. In this case, set
\[
H_3(q)=\sup \bigl\{s \le H_2(q)\dvtx|X_s-q|
\le\tfrac{1}{2}t^{({\alpha-d+2})/({4\alpha})} (\log t)^{{1}/{2}+\delta} \bigr\}.
\]
Then for $s \in[H_3(q),H_2(q)]$ and $c<C(d,\alpha)/4$,
\begin{eqnarray*}
V_{\omega}(X_s) &\ge& V_{\omega}(m_{L_t})+
ct^{-({\alpha-d+2})/({2\alpha})}(\log t)^{1+2\delta}
\\
&\ge&\lambda_1^{\omega}\bigl((-t,t)^d\bigr)+
\frac{c}{2}t^{-({\alpha-d+2})/({2\alpha})}(\log t)^{1+2\delta}
\end{eqnarray*}
by Proposition~\ref{PC-L_t} and (\ref{ev-bound}).
Let us define $E_{k,l}(q)=\{H_3(q) \in[k,k+1),\break H_2(q)\in[l,l+1)\}$ for
integers $0 \le k \le l \le t-1$.
Suppose first that $l-k> t^{({\alpha-d+2})/({2\alpha})}$ (slow crossing).
Then using the Markov property at time $l$ and then (\ref{unif-FK})
exactly the same way as in the proof of
Proposition~\ref{weak-localization}, we obtain
%
%
\begin{eqnarray}
\label{slow-crossing-2} && \E\otimes E_0 \biggl[\exp \biggl\{-\int
_0^{t} V_{\omega}(X_s)\,\d s
\biggr\}\dvtx E_{k,l}(q) \biggr]
\nonumber
\\
&&\qquad \le\E \biggl[c(d)^2 \bigl(1+\bigl(t\lambda_1^{\omega}
\bigl((-t,t)^d\bigr)\bigr)^{{d}/{2}} \bigr)^2 \exp
\bigl\{-(k+1)\lambda_1^{\omega}\bigl((-t,t)^d\bigr)
\bigr\}
\nonumber
\\
&&\hspace*{10.7pt}\qquad\quad{}\times\exp \biggl\{-(l-k-1) \biggl(\lambda_1^{\omega}
\bigl((-t,t)^d\bigr)+\frac{c}{2}t^{-({\alpha-d+2})/({2\alpha})} (\log
t)^{1+2\delta} \biggr) \biggr\}
\\
&&\hspace*{179.3pt}\qquad\quad{}\times\exp \bigl\{-(t-l)\lambda_1^{\omega}\bigl((-t,t)^d
\bigr) \bigr\} \biggr]
\nonumber
\\
&&\qquad =\O \bigl(\exp \bigl\{-(\log t)^{1+2\delta} \bigr\} \bigr) \E \bigl[
\exp \bigl\{-t\lambda_1^{\omega}\bigl((-t,t)^d\bigr)
\bigr\} \bigr].\nonumber
\end{eqnarray}
For the case $l-k\le t^{({\alpha-d+2})/({2\alpha})}$ (fast crossing),
we again proceed as in the proof of (\ref{crossing}) to see
%
%
\begin{eqnarray}
\label{fast-crossing-2} && \E \biggl[\exp \biggl\{-\int_0^{t}
V_{\omega}(X_s)\,\d s \biggr\}\dvtx E_{k,l}(q)
\biggr]
\nonumber
\\
&&\qquad \le\E \bigl[c(d)^2 \bigl(1+\bigl(t\lambda_1^{\omega}
\bigl((-t,t)^d\bigr)\bigr)^{{d}/{2}} \bigr)^2 \exp
\bigl\{-t\lambda_1^{\omega}\bigl((-t,t)^d\bigr)\bigr
\}
\nonumber\\[-8pt]\\[-8pt]
&&\qquad\quad\hspace*{7.3pt}{}\times P_0 \bigl( \sup \bigl\{|X_s|\dvtx0\le s \le
t^{({\alpha-d+2})/({2\alpha})} \bigr\} >t^{({\alpha-d+2})/({4\alpha})} (\log t)^{{1}/{2}+\delta} \bigr) \bigr]
\nonumber
\\
&&\qquad \le\O \bigl(\exp \bigl\{-(\log t)^{1+2\delta} \bigr\} \bigr) \E
\bigl[\exp \bigl\{-t\lambda_1^{\omega}\bigl((-t,t)^d
\bigr)\bigr\} \bigr].\nonumber
\end{eqnarray}
Summing (\ref{slow-crossing-2}) and (\ref{fast-crossing-2}) over
$0 \le k \le l \le t-1$ and using Lemma~\ref{lem5},
we find that
\[
Q_t \bigl(H_1(q) < H_2(q) \le t,
m_{L_t}\in q+[0,1 )^d \bigr) =\O \bigl(\exp \bigl\{-(\log
t)^{1+2\delta} \bigr\} \bigr).
\]
The other case $H_2(q)<H_1(q)\le t$ can also
be treated in the same way as above by using
\begin{eqnarray*}
\tilde H_3(q)&=&\sup \bigl\{s \in\bigl[0, H_1(q)\bigr)\dvtx
|X_s-q|> t^{({\alpha-d+2})/({4\alpha})}(\log t)^{{1}/{2}+\delta} \bigr\},
\\
\tilde H_2(q)&=&\inf \bigl\{s \ge\tilde H_3(q)\dvtx
|X_s-q|\le\tfrac
{1}{2}t^{({\alpha-d+2})/({4\alpha})}(\log t)^{{1}/{2}+\delta
}
\bigr\}
\end{eqnarray*}
instead of $H_3(q)$ and $H_2(q)$, respectively.
Consequently, we obtain
\[
Q_t \bigl(H_2(q) \le t, m_{L_t}\in q+[0,1
)^d \bigr) =\O \bigl(\exp \bigl\{-(\log t)^{1+2\delta} \bigr\}
\bigr).
\]

Since the possible values of $q$ are only polynomially
many due to the weak localization result,
we can sum over $q$ to see
\[
\lim_{t\to\infty}Q_t \Bigl(\sup_{s\in[0,t]}
|X_s-m_{L_t}|\le2t^{({\alpha-d+2})/({4\alpha})}(\log t)^{
{1}/{2}+\delta}
\Bigr) = 1.
\]
Finally observe that on the event in the left-hand side, we have
$|m_{L_t}|\le2t^{({\alpha-d+2})/({4\alpha})}(\log t)^{
{1}/{2}+\delta}$
since $X_0=0$. Therefore we arrive at the desired conclusion,
\[
\lim_{t\to\infty} Q_t \Bigl(\sup_{s\in[0,t]}|X_s|
\le4t^{({\alpha-d+2})/({4\alpha})} (\log t)^{{1}/{2}+\delta} \Bigr)=1.
\]
\upqed\end{pf*}
%
By the strong localization result, it, in particular, follows that
$m_{L_t}$ is close to $m_t(\omega)$, which completes the proof
of Theorem~\ref{potential-confinement}.
\begin{pf*}{Proof of Theorem~\ref{potential-confinement}}
Fix $\varepsilon>0$ so small that Proposition~\ref{PC-L_t} holds.
Then, for sufficiently large $t$,
\[
t^{({\alpha-d+2})/({4\alpha})-{\varepsilon}/{4}} <|x-m_{L_t}|<2t^{({\alpha-d+2})/({4\alpha})}\log t
\]
implies $V_{\omega}(x)>V_{\omega}(m_{L_t})$. Therefore, $m_{L_t}$ is
within $t^{({\alpha-d+2})/({4\alpha})-{\varepsilon}/{4}}$ from the
minimizer of $V_{\omega}$ in $B(m_{L_t},
2t^{({\alpha-d+2})/({4\alpha})}\log t)$. On the other\vspace*{1pt}
hand, Theorem~\ref{localization} implies $m_{L_t}\in
B(0,t^{({\alpha-d+2})/({4\alpha})}\log t)$, and thus the above
minimizer is nothing but $m_t(\omega)$. It is now easy to deduce
Theorem~\ref{potential-confinement} from Proposition~\ref{PC-L_t}.
\end{pf*}
%
\section{Scaling limit of the occupation time measure}
\label{OT}
We prove Theorem~\ref{occupation-time} in this section.
Given Theorems~\ref{localization} and~\ref{potential-confinement},
it is more or less straightforward. Indeed, what we do is,
replacing $V_{\omega}$ by a quadratic function, using the
Girsanov formula, and applying the large deviation principle
for the Ornstein--Uhlenbeck process.
\begin{pf*}{Proof of Theorem~\ref{occupation-time}}
Let us introduce the event
\[
G_1= \Bigl\{\sup_{0 \le u \le t}|X_u| < r(t) (\log
t)^{{3}/{4}} \Bigr\}.
\]
By Theorem~\ref{localization} and Proposition~\ref{PC-L_t},
we know that
\begin{eqnarray*}
&& Q_t (\Tilde{L}_t \in A )
\\
&&\qquad =Q_t \bigl(\Tilde{L}_t \in A,
G_1\cap\mathtt{PC}(L_t) \bigr)+o(1)
\\
&&\qquad =\frac{1}{Z_t}\E\otimes E_0 \biggl[\exp \biggl\{-
\int_0^{t} V_{\omega}(X_s)\,
\d s \biggr\}\dvtx\Tilde{L}_t \in A, G_1\cap\mathtt{PC}(L_t) \biggr]+o(1),
\end{eqnarray*}
and thus we may concentrate on $G_1\cap\mathtt{PC}(L_t)$.
On this event, we have
\begin{eqnarray*}
&&\exp \biggl\{-\int_0^{t}
V_{\omega}(X_s)\,\d s \biggr\}1_{\mathtt{PC}(L_t)\cap G_1}
\\
&&\qquad \le\exp \bigl\{-tV_{\omega}(m_{L_t}) +o
\bigl(t^{({\alpha+d-2})/({2\alpha})} \bigr) \bigr\}
\\
&&\qquad\quad\hspace*{0pt}{}\times\exp \biggl\{-C(d,\alpha) \int_0^{tr(t)^{-2}} |
\Tilde{X}_s-m_{\Tilde L_t}|^2 \,\d s \biggr
\}1_{G_1}.
\end{eqnarray*}
Let $q\in t^{-1}\Z^d$ and $f\in C_c(\R^d)$.
We integrate the above over
\[
F(q,f)= \biggl\{\|m_{\Tilde L_t}-q\|_{\infty}\le(2t)^{-1},
\biggl\llvert\int f \,\d\Tilde{L}_t-\int f \,\d\nu_q
\biggr\rrvert>\varepsilon \biggr\}.
\]
Note that if $|q|>2\log t$, then $F(q,f)\cap G_1=\varnothing$.
Hence we assume $|q|\le2\log t$ in what follows.
Since
\[
\int_0^{tr(t)^{-2}}|\Tilde{X}_s-m_{\Tilde L_t}|^2
\,\d s =\int_0^{tr(t)^{-2}}|\Tilde{X}_s-q|^2
\,\d s +o \bigl(t^{({\alpha+d-2})/({2\alpha})} \bigr)
\]
on $G_1 \cap F(q,f)$, we have
\begin{eqnarray*}
&&\E\otimes E_0 \biggl[\exp \biggl\{-\int_0^{t}
V_{\omega}(X_s)\,\d s \biggr\}\dvtx\mathtt{PC}(L_t)
\cap G_1 \cap F(q,f) \biggr]
\\
&&\qquad \le\exp \bigl\{-a_1t^{{d}/{\alpha}} +o
\bigl(t^{({\alpha+d-2})/({2\alpha})} \bigr) \bigr\}
\\
&&\qquad\quad\hspace*{0pt}{}\times E_0 \biggl[\exp \biggl\{-C(d,\alpha) \int_0^{tr(t)^{-2}}
|\Tilde{X}_s-q|^2 \,\d s \biggr\}\dvtx G_1
\cap F(q,f) \biggr].
\end{eqnarray*}
By the scaling invariance of Brownian motion and the Girsanov formula,
it follows that the right-hand side equals
\[
\exp \bigl\{-a_1t^{{d}/{\alpha}} - \bigl(a_2+o(1)
\bigr)t^{({\alpha+d-2})/({2\alpha})} \bigr\} R_0^q \biggl[
\frac{\psi(-q)} {
\psi(X_{tr(t)^{-2}}-q)}\dvtx\tilde G_1 \cap\tilde F(q,f) \biggr],
\]
where
\begin{eqnarray*}
\psi(x)&=&\exp \Biggl\{-\sqrt{\frac{C(d,\alpha)}{2}}|x|^2 \Biggr\},
\\
\tilde G_1&=& \Bigl\{\sup_{s\in[0,tr(t)^{-2}]}|X_s| \le\log
t \Bigr\},
\\
\tilde F(q,f)&=& \biggl\{\|m_{L_{tr(t)^{-2}}}-q\|_{\infty}
\le(2t)^{-1}, \biggl\llvert\int f \,\d L_{tr(t)^{-2}}-\int f \,\d
\nu_q \biggr\rrvert>\varepsilon \biggr\}.
\end{eqnarray*}
Recalling $|q|\le2\log t$, we have
\[
\frac{\psi(-q)}{\psi(X_{tr(t)^{-2}}-q)} =\O \bigl(\exp \bigl\{(\log t)^2 \bigr\} \bigr)
\]
on $\tilde G_1$.
It is well known that $L_s$ under $R_0^q$ satisfies the full
large deviation principle with rate $s$ in the space of probability
measures equipped with weak topology,
and the rate function has unique zero at $\nu_q$;
see, for example, \cite{DS89}.
Therefore, we arrive at
\begin{eqnarray*}
&&\frac{1}{Z_t} \E\otimes E_0 \biggl[\exp \biggl\{-\int
_0^{t} V_{\omega}(X_s)\,\d s
\biggr\}\dvtx\mathtt{PC}(L_t)\cap G_1 \cap F(q,f) \biggr]
\\
&&\qquad \le\exp \bigl\{o \bigl(t^{({\alpha+d-2})/({2\alpha})} \bigr) \bigr\}
R_0^q \biggl( \biggl\llvert\int f \,\d
L_{tr(t)^{-2}}-\int f \,\d\nu_q \biggr\rrvert>\varepsilon \biggr)
\\
&&\qquad \le\exp \bigl\{\delta(\varepsilon)t^{({\alpha+d-2})/({2\alpha})} \bigr\}
\end{eqnarray*}
for some $\delta(\varepsilon)>0$.
Summing over $q \in t^{-1}\Z^d \cap B(0, 2\log t)$ and replacing
$m_{L_t}$ by $m_t(\omega)$ using Theorem~\ref{potential-confinement},
we obtain (\ref{thm4-1}).
\end{pf*}
%
\section{Fluctuation of the local minimum of the potential}
\label{fluctuation}
We prove Theorem~\ref{local-minimum} in this section.
The argument\vspace*{1pt} is similar to that in Section~\ref{weak-to-PC}, but
we need a better control on
$\int|x-m_{\Tilde{L}_t}|^2 \Tilde{L}_t(\d x)$.
We deduce it from Theorem~\ref{occupation-time}.
\begin{pf*}{Proof of Theorem~\ref{local-minimum}}
We first prove
%
%
\begin{equation}
\label{thm4-2} \lim_{t \to\infty}Q_t \biggl( \biggl\llvert
\int|x-m_{\Tilde{L}_t}|^2 \Tilde{L}_t(\d x)-
\bigl(8C(d,\alpha) \bigr)^{-{1}/{2}} \biggr\rrvert>\varepsilon \biggr)=0
\end{equation}
for any $\varepsilon>0$.
To this end, recall that
Corollary~\ref{variance-bound} allows us to assume
\[
\int|x-m_{\Tilde{L}_t}|^2 \Tilde{L}_t(\d x)\le
M_2.
\]
Under this assumption,
\[
\biggl\llvert\int|x-m_{\Tilde{L}_t}|^2 \Tilde{L}_t(
\d x)- \int|x-m_{\Tilde{L}_t}|^2 \nu_{m_{\Tilde{L}_t}}(\d x) \biggr
\rrvert>\varepsilon
\]
implies
\[
\biggl\llvert\int \bigl(|x-m_{\Tilde{L}_t}|^2\wedge M \bigr)
\Tilde{L}_t(\d x)- \int \bigl(|x-m_{\Tilde{L}_t}|^2
\wedge M \bigr) \nu_{m_{\Tilde{L}_t}}(\d x) \biggr\rrvert>\frac
{\varepsilon}{2}
\]
for sufficiently large $M$.
Since $\int|x|^2 \nu_0(\d x)=(8C(d,\alpha))^{-{1}/{2}}$,
(\ref{thm4-2}) follows from Theorem~\ref{occupation-time}.
Now let us define the class of probability measures on $\R^d$ by
\begin{eqnarray*}
&&
\mathcal{P}^{\varepsilon}_{\mathrm{loc}}:= \biggl\{\mu\dvtx
\operatorname{supp}\mu \subset B \bigl(0,t^{({\alpha-d+2})/({4\alpha})}\log t \bigr),
\\
&&\hspace*{41pt}\biggl\llvert\int|x-m_{\mu}|^2\mu(\d x)- \bigl(8C(d,\alpha)
\bigr)^{-{1}/{2}} \biggr\rrvert<\varepsilon t^{-({\alpha-d+2})/({2\alpha})} \biggr\}.
\end{eqnarray*}
By Theorem~\ref{localization} and (\ref{thm4-2}), we know that
$\lim_{t\to\infty}Q_t(L_t\notin\mathcal{P}^{\varepsilon}_{\mathrm{loc}})=0$
for any $\varepsilon>0$.

Now let us prove Theorem~\ref{local-minimum}(i).
It suffices to show the assertion with $m_t(\omega)$ replaced
by $m_{L_t}$ in view of (\ref{bottom-location})
and Lemma~\ref{lem7}.
For $\mu\in\mathcal{P}^{\varepsilon}_{\mathrm{loc}}$,
the same argument as for Lemma~\ref{lem7} yields
%
%
\begin{eqnarray}
\label{weighted-mean} &&\E_t^{\mu} \bigl[V_{\omega}(m_{\mu})
\bigr]
\nonumber
\\
&&\qquad =\int\hat v(m_{\mu}-y) e^{-t\int\hat v(z-y)\mu(\d z)}\,\d y
\nonumber\\[-8pt]\\[-8pt]
&&\qquad =a_1\frac{d}{\alpha}t^{-({\alpha-d})/{\alpha}}\nonumber\\
&&\qquad\quad{} - \biggl(
\frac{\alpha-d+2}{\alpha}C(d,\alpha)+o(1) \biggr) t^{-({\alpha-d+2})/{\alpha}} \int|x-m_{\mu}|^2
\mu(\d x),\nonumber
\end{eqnarray}
where $o(1)$ is uniform in $\mu$.
Then it follows
\begin{eqnarray*}
&& Q_t \bigl[V_{\omega}(m_{L_t}) \bigr]
\\
&&\qquad\sim\frac{1}{Z_t}E_0 \bigl[\E \bigl[e^{-t\langle L_t,V_{\omega}\rangle
}
\bigr] \E_t^{L_t} \bigl[V_{\omega}(m_{L_t})
\bigr] \dvtx L_t\in\mathcal{P}^{\varepsilon}_{\mathrm{loc}} \bigr]
\\
&&\qquad=a_1\frac{d}{\alpha}t^{-({\alpha-d})/{\alpha}} -\frac{\alpha
-d+2}{\alpha}
\sqrt{ \frac{C(d,\alpha)}{8}}t^{-({\alpha-d+2})/({2\alpha})}\\
&&\qquad\quad{} +\varepsilon O \bigl(t^{-({\alpha-d+2})/({2\alpha})}
\bigr),
\end{eqnarray*}
and letting $\varepsilon\downarrow0$, we get
Theorem~\ref{local-minimum}(i).

Let us turn to Theorem~\ref{local-minimum}(ii).
It suffices to consider the assertion with $m_t(\omega)$ replaced
by $m_{L_t}$ again. Indeed,
as is mentioned in Remark~\ref{rem1}, we have
\[
\lim_{t\to\infty} Q_t \bigl(\bigl|V_{\omega}
\bigl(m_t( \omega) \bigr)-V_{\omega}(m_{L_t}) -
\E_t^{\mu} \bigl[V_{\omega} \bigl(m_t(
\omega) \bigr)-V_{\omega}(m_{L_t}) \bigr]\bigr| <t^{-\beta}
\bigr)= 1
\]
for any
\[
\frac{2\alpha-d}{2\alpha}<\beta<\min \biggl\{\frac{3\alpha-d+2}{4\alpha
},\frac{\alpha
-d+6}{4\alpha}
\biggr\}.
\]
Now it follows by a well-known formula for the characteristic functional
for Poisson point process that for
$\mu\in\mathcal{P}^{\varepsilon}_{\mathrm{loc}}$,
\begin{eqnarray*}
&&\log\E_t^{\mu} \bigl[\exp \bigl\{i\theta
t^{({2\alpha
-d})/({2\alpha})} \bigl( V_{\omega}(m_{\mu})-\E_t^{\mu}
\bigl[V_{\omega
}(m_{\mu}) \bigr] \bigr) \bigr\} \bigr]
\\
&&\qquad =\int \bigl(e^{
i\theta t^{({2\alpha-d})/({2\alpha})}
\hat v(m_{\mu}-y)}-1-i\theta t^{({2\alpha-d})/({2\alpha})} \hat
v(m_{\mu}-y) \bigr)\\
&&\qquad\quad\hspace*{8.2pt}{}\times e^{-t\int\hat v(z-y)\mu(\d z)}\,\d y.
\end{eqnarray*}
Using Taylor's theorem, replacing $\int\hat v(z-y)\mu(\d z)$
by $v(m_{\mu}-y)$, and changing variable, one can show that
\begin{eqnarray*}
&& \log\E_t^{\mu} \bigl[\exp \bigl\{i\theta
t^{({2\alpha
-d})/({2\alpha})} \bigl( V_{\omega}(m_{\mu})-\E_t^{\mu}
\bigl[V_{\omega
}(m_{\mu}) \bigr] \bigr) \bigr\} \bigr]
\\
&&\qquad \sim-\frac{\theta^2}{2}\int t^{({2\alpha-d})/{\alpha}} |y|^{-2\alpha}
e^{-t|y|^{-\alpha}}\,\d y
\\
&&\qquad =-\frac{\theta^2}{2} \alpha\sigma_d\Gamma \biggl(
\frac{3\alpha-d+1}{\alpha} \biggr).
\end{eqnarray*}
We leave the details to the reader.
From this, Theorem\ref{local-minimum}(ii) follows by the same
way as above.
\end{pf*}
%
\section{Scaling limit of the process}
\label{SL}
In this section, we prove Theorem~\ref{scaling-limit}.
Given the strong localization and the potential confinement,
we can basically follow the argument in \cite{Szn91b}.
We write $B_t$ for $B(m_t(\omega),r(t)\log t)$.
\begin{pf*}{Proof of Theorem~\ref{scaling-limit}}
We define the good events by
\begin{eqnarray*}
G_1(s) &=& \Bigl\{\sup_{0 \le u \le s}|X_u| < r(t) (
\log t)^{{3}/{4}} \Bigr\},
\\
G_2 &=& \Bigl\{\sup_{x \in B_t}\bigl|V_{\omega}(x)-V_{\omega}
\bigl(m_t(\omega) \bigr) -p_t \bigl(x-m_t(
\omega) \bigr)\bigr| \le t^{-({\alpha-d+2})/({2\alpha})-\varepsilon_0},
\\
&&\hspace*{186.1pt} \bigl|m_t(\omega)\bigr| \le r(t) (\log t)^{{3}/{4}} \Bigr\}.
\end{eqnarray*}
Note that $X_{[0,t]}\subset B_t$ on $G_1(t)\cap G_2$.
Due to Theorems~\ref{localization} and~\ref{potential-confinement},
it suffices to show that
for any $T>0$ and $f$, a bounded continuous function on
$C([0,T], \R^d)$,
\begin{eqnarray*}
&&\lim_{t\to\infty}Q_t \bigl[f(\Tilde{X})\dvtx
G_1(t)\cap G_2 \bigr]
\\
&&\qquad = \int\d m \biggl(\frac{\sqrt{C(d,\alpha)}}{\sqrt{2}\pi} \biggr)^{{d}/{2}} \exp
\Biggl\{- \sqrt{\frac{C(d,\alpha)}{2}}|m|^2 \Biggr\} R^m_0
\bigl[f(X) \bigr]
\end{eqnarray*}
with an obvious abuse of notation.
Let $\lambda_i^{\omega}(B_t)$ denote the $i$th smallest eigenvalue
of $-\Delta/2+V_{\omega}(x)-V_{\omega}(m_t({\omega}))$ in $B_t$
with the Dirichlet boundary condition and $\phi_{\omega}$
the $L^2$-normalized nonnegative principal eigenfunction.
Let us further denote by $\psi_{\omega}$ the
$L^2$-normalized nonnegative principal eigenfunctions of
$-\Delta/2+p_t(x-m_t({\omega}))$ on $\R^d$.
%
%
\begin{lemma}
\label{lem9}
\textup{(i)}
$\lambda_i^{\omega}(B_t)
=V_{\omega}(m_t(\omega))+(a_2+O(t^{-\varepsilon_0}))
t^{-({\alpha-d+2})/({2\alpha})}$
uniformly in $i \in\{1,2\}$ and $\omega\in G_2$.

\textup{\hphantom{i}(ii)} There exists $c_3>0$ such that for any $\omega\in G_2$,
\[
\lambda_2^{\omega}(B_t)-\lambda^1_{\omega}(B_t)
\ge c_3 t^{-({\alpha-d+2})/({2\alpha})}.
\]

\textup{(iii)} ${\sup_{\omega\in G_2}}\|\phi_{\omega}-\psi_{\omega}\|_2 \to0$
as $t \to\infty$.
\end{lemma}
%
%
\begin{pf}
It may be assumed that $m_t(\omega)=0$ by a spatial shift.
Let $\lambda_i^{p_t}(U)$ denote the
$i$th smallest eigenvalue of $-\Delta/2+p_t(x)$ in $U$
with the Dirichlet boundary condition.
Recall that when $U=\R^d$, the corresponding eigenfunctions
are products of the Hermite polynomials multiplied by the
Gaussian density with variance $\sqrt{r(t)/2C(d,\alpha)}$.
In particular, they are of\break $\O(\exp\{-(\log t)^2\})$
near $\partial B_t$.
Thus, by using the Rayleigh--Ritz variational formula
and a standard cut-off argument, one can show that
%
%
\begin{equation}
\label{ev-approx1} \lambda_i^{p_t}(B_t) =
\lambda_i^{p_t} \bigl(\R^d \bigr) +\O \bigl(
\exp \bigl\{-(\log t)^2 \bigr\} \bigr).
\end{equation}
Since
%
%
\begin{equation}
\label{ev-approx2} \bigl\llvert \bigl(\lambda_i^{\omega}(B_t)-V_{\omega}
\bigl(m_t(\omega) \bigr) \bigr)- \lambda_i^{p_t}(B_t)
\bigr\rrvert\le t^{-({\alpha-d+2})/({2\alpha})-\varepsilon_0}
\end{equation}
on $G_2$, the first assertion follows.
Moreover, we know
%
%
\begin{eqnarray}
\label{spectral-gap} \lambda_2^{p_t} \bigl(
\R^d \bigr)-\lambda_1^{p_t} \bigl(
\R^d \bigr) &=& r(t)^{-2}\lambda_2^{p_1}
\bigl(\R^d \bigr) -r(t)^{-2}\lambda_1^{p_1}
\bigl( \R^d \bigr)
\nonumber\\[-8pt]\\[-8pt]
&\ge& c_3 t^{-({\alpha-d+2})/({2\alpha})}\nonumber
\end{eqnarray}
for some $c_3>0$ by a scaling and the fact that
$-\Delta/2+C(d,\alpha)|x|^2$ has positive spectral gap.
From (\ref{ev-approx1})--(\ref{spectral-gap}), the second
assertion follows.

Next, we introduce the Dirichlet form
$(\mathcal{E}_{\omega}, \mathcal{D}(\mathcal{E}_{\omega}))$
associated with $-\Delta/2+p_t(x)$ on $L^2(\R^d)$.
Then $\phi_{\omega} \in\mathcal{D}(\mathcal{E}_{\omega})$
since $\phi_{\omega}$ has compact support and\break
$\int|\nabla\phi_{\omega}|^2\,\d x<\infty$.
Let us decompose $\phi_{\omega}$ as
$\gamma_t\psi_{\omega}+\sqrt{1-\gamma_t^2}\psi_2$
by using $\psi_2 \in\mathcal{D}(\mathcal{E}_{\omega})$
with unit $L^2$-norm and
orthogonal to $\psi_{\omega}$ with respect to
$\mathcal{E}_{\omega}(\cdot,\cdot)+(\cdot,\cdot)_{L^2}$.
Note that we have
\[
\mathcal{E}_{\omega}(\psi_{\omega},\psi_2) = \bigl((-
\Delta/2+p_t)\psi_{\omega},\psi_2
\bigr)_{L^2} =\lambda_1^{p_t} \bigl(
\R^d \bigr) (\psi_{\omega},\psi_2)_{L^2}
\]
since $\psi_{\omega}$ is an eigenfunction.
Hence in fact $\mathcal{E}_{\omega}(\psi_{\omega},\psi_2)=
(\psi_{\omega},\psi_2)_{L^2}=0$, and it follows that
%
%
\begin{eqnarray}
\label{form} \mathcal{E}_{\omega}(\phi_{\omega},\phi_{\omega})
&=& \gamma_t^2\mathcal{E}_{\omega}(
\psi_{\omega}, \psi_{\omega}) + \bigl(1-\gamma_t^2
\bigr) \mathcal{E}_{\omega}(\psi_2,\psi_2)
\nonumber\\[-8pt]\\[-8pt]
&=& \gamma_t^2 \lambda_1^{p_t}
\bigl(\R^d \bigr) + \bigl(1-\gamma_t^2
\bigr) \mathcal{E}_{\omega}(\psi_2,\psi_2).\nonumber
\end{eqnarray}
Now for $\omega\in G_2$, one can easily find
\begin{eqnarray*}
&& \mathcal{E}_{\omega}(\phi_{\omega},\phi_{\omega})
\\
&&\qquad=\frac{1}{2}\int\bigl|\nabla\phi_{\omega}(x)\bigr|^2\,\d x\\
&&\qquad\quad{}
+\int \bigl(V_{\omega}(x)-V_{\omega}(0) \bigr)\phi_{\omega}(x)^2
\,\d x +o \bigl(t^{-({\alpha-d+2})/({2\alpha})} \bigr)
\\
&&\qquad=\lambda_1^{\omega}(B_t)
-V_{\omega}(0) +o \bigl(t^{-({\alpha-d+2})/({2\alpha})} \bigr)
\\
&&\qquad= \lambda_1^{p_t} \bigl(\R^d \bigr) +o
\bigl(t^{-({\alpha-d+2})/({2\alpha})} \bigr),
\end{eqnarray*}
where the last line is due to part (i).
On the other hand, by the variational formula and (\ref{spectral-gap}),
\[
\mathcal{E}_{\omega}(\psi_2,\psi_2) \ge
\lambda_2^{p_t} \bigl(\R^d \bigr) \ge
\lambda_1^{p_t} \bigl(\R^d \bigr) +
c_3 t^{-({\alpha-d+2})/({2\alpha})}.
\]
Substituting these relations into (\ref{form}), we obtain
$\gamma_t^2 \to1$ as $t \to\infty$.
Finally, since both $\phi_{\omega}$ and $\psi_{\omega}$ are
nonnegative, $\gamma_t$ must converge to 1, and the last
assertion is proved.
\end{pf}
Let us define
$\tau=t^{({\alpha-d+2})/({2\alpha})+{\varepsilon_0}/{2}}$
and
\begin{eqnarray*}
F_1(x,\omega)&=&E_x \biggl[\exp \biggl\{-\int
_0^{t-\tau} V_{\omega
}(X_s)\,\d
s \biggr\}\dvtx G_1(t-\tau) \biggr],
\\
F_2(x,\omega)&=&\langle\phi_{\omega}, 1 \rangle
\phi_{\omega}(x) \exp \bigl\{-(t-\tau)\lambda_1^{\omega}(B_t)
\bigr\}
\end{eqnarray*}
for the ease of notation.
Note that $\tau>Tr(t)^2$ for sufficiently large $t$.
Applying the Markov property at time $\tau$, we find
that the numerator of $Q_t[f(\Tilde{X})\dvtx\break G_1(t)\cap G_2]$ is
%
%
\begin{eqnarray}
\label{to-replace} && \E\otimes E_0 \biggl[f(\Tilde{X})\exp \biggl
\{- \int_0^{t} V_{\omega
}(X_s)
\,\d s \biggr\}\dvtx G_1(t)\cap G_2 \biggr]
\nonumber
\\
&&\qquad =\E \biggl[E_0 \biggl[f(\Tilde{X})\exp \biggl\{-\int
_0^{\tau
} V_{\omega}(X_s)\,\d s
\biggr\}
\\
&&\qquad\quad\hspace*{65.3pt}{}\times F_1(X_{\tau},\omega)\dvtx G_1(\tau) \biggr]
\dvtx G_2 \biggr].\nonumber
\end{eqnarray}
To replace $F_1$ by $F_2$ in this expression,
we estimate the difference as
%
%
\begin{eqnarray}
\label{projection} && \biggl| E_0 \biggl[f(\Tilde{X})\exp \biggl\{-\int
_0^{\tau} V_{\omega}(X_s)\,\d s
\biggr\}
\nonumber
\\
&&\quad\hspace*{7.7pt}{}\times \bigl(F_1(X_{\tau},\omega)-F_2(X_{\tau},
\omega) \bigr)\dvtx G_1(\tau) \biggr] \biggr|
\nonumber\\[-8pt]\\[-8pt]
&&\hspace*{7.7pt}\qquad \le\|f\|_{\infty} \int p(1,0,x)S^{\omega}_{\tau-1}
\bigl|F_1(\cdot,\omega)-F_2(\cdot,\omega)\bigr|(x)\,\d x
\nonumber
\\
&&\hspace*{7.7pt}\qquad \le\|f\|_{\infty} \biggl(\int p(1,0,x)^2 \,\d x
\biggr)^{{1}/{2}} \bigl\|S^{\omega}_{\tau-1}\bigl|F_1(\cdot,
\omega)-F_2(\cdot,\omega)\bigr|\bigr\|_2,\nonumber
\end{eqnarray}
where $p(t,x,y)$ denotes the transition kernel of the Brownian
motion and $\{S^{\omega}_s\}_{s \ge0}$ the semigroup
generated by $-\Delta/2+V_{\omega}$ in $B_t$ with
the Dirichlet boundary condition.
Now, it is well known that
$\|S^{\omega}_s\|_{L^2\to L^2}=\exp\{-s\lambda_1^{\omega}(B_t)\}$
and also
by considering eigenfunction expansion, one can deduce
from Lem\-ma~\ref{lem9}(ii) that
\begin{eqnarray*}
&&\bigl\llVert F_1(\cdot,\omega)-F_2(\cdot,\omega) \bigr
\rrVert_2 \\
&&\qquad\le|B_t|^{{1}/{2}} \exp \bigl\{-(t-\tau)
\bigl(\lambda_1^{\omega}(B_t) +c_3t^{-({\alpha-d+2})/({2\alpha})}
\bigr) \bigr\}
\end{eqnarray*}
on $G_2$.
Therefore for sufficiently large $t$,
%
%
\begin{equation}
\label{F_1-F_2} \mbox{RHS of (\ref{projection})} \le\exp
\biggl\{-t\lambda_1^{\omega}(B_t) -
\frac{c_3}{2}t^{({\alpha+d-2})/({2\alpha})} \biggr\}.
\end{equation}
Coming back to (\ref{to-replace}) and replacing
$F_1$ by $F_2$, we arrive at
%
%
\begin{eqnarray}
\label{replaced} \mbox{LHS of (\ref{to-replace})} &=&\E \biggl[E_0
\biggl[f(\Tilde{X})\phi_{\omega}(X_{\tau})\exp \biggl\{-\int
_0^{\tau} V_{\omega}(X_s)\,\d s
\biggr\}\dvtx G_1(\tau) \biggr]
\nonumber
\\
&&\hspace*{110.5pt}{}\times\langle\phi_{\omega}, 1\rangle e^{-(t-\tau)\lambda_1^{\omega}(B_t)} \dvtx G_2
\biggr]
\\
&&{}+\O \bigl(\exp \bigl\{-t^{({\alpha+d-2})/({2\alpha})} \bigr\} \bigr) \E \bigl[ \exp \bigl
\{-t\lambda_1^{\omega}(B_t) \bigr\} \bigr].\nonumber
\end{eqnarray}
Recall that the last term is of $o(Z_t)$ due to Lemma~\ref{lem5}.
%
%
\begin{lemma}
\label{lem10}
Uniformly in $\omega\in G_2$,
%
%
\begin{eqnarray}
\label{lemma10} && E_0 \biggl[f(\Tilde{X})\phi_{\omega}(X_{\tau})
\exp \biggl\{-\int_0^{\tau} V_{\omega}(X_s)
\,\d s \biggr\}\dvtx G_1(\tau) \biggr]
\nonumber\\[-8pt]\\[-8pt]
&&\qquad =e^{-\tau\lambda_1^{\omega}(B_t)+o(1)} \bigl( \psi_{\omega
}(0)R_0^{\tilde{m}_t(\omega)}
\bigl[f(X) \bigr] +\O \bigl(\exp \bigl\{-(\log t)^2 \bigr\} \bigr)
\bigr).\nonumber
\end{eqnarray}
\end{lemma}
%
%
\begin{pf}
Let us use the orthogonal decomposition
$\phi_{\omega}=\gamma_t\psi_{\omega}+\sqrt{1-\gamma_t^2}\psi_2$
in the proof of Lemma~\ref{lem9}.
Denoting the left-hand side of (\ref{lemma10}) by $T(\phi_{\omega})$,
we have
%
%
\begin{equation}
\label{decomposition} T(\phi_{\omega})=\gamma_tT(
\psi_{\omega}) +\sqrt {1-\gamma_t^2}T(
\psi_2).
\end{equation}
We begin with the first term in the right-hand side.
We know $\gamma_t \to1$ by Lem\-ma~\ref{lem9}(iii)
and on $G_2$,
\begin{eqnarray*}
T(\psi_{\omega}) &=& E_0 \biggl[f(\Tilde{X})
\psi_{\omega}(X_{\tau}) \exp \biggl\{-\int_0^{\tau}
V_{\omega}(X_s)\,\d s \biggr\}\dvtx G_1(\tau)
\biggr]
\\
&\le& \exp \bigl\{-\tau V_{\omega} \bigl(m_t(\omega)
\bigr)+t^{-{\varepsilon
_0}/{2}} \bigr\}
\\
&&{}\times E_0 \biggl[f(\Tilde{X})\psi_{\omega}(X_{\tau})
\exp \biggl\{-\int_0^{\tau} p_t
\bigl(X_s-m_t(\omega) \bigr)\,\d s \biggr\} \dvtx
G_1(\tau) \biggr]
\end{eqnarray*}
thanks to our choice of $\tau$.
On the other hand, it follows by the Girsanov transform
and Brownian scaling that
\begin{eqnarray*}
&& R_0^{\tilde{m}_t(\omega)} \bigl[f(X)\dvtx X_{[0, t^{\varepsilon_0}]} \subset
r(t)^{-1}B_t \bigr]
\\
&&\qquad =e^{a_2} E_0 \biggl[f(X)\frac{\psi_{\omega}(r(t)X_{r(t)^{-2}\tau})} {
\psi_{\omega}(0)}
\\
&&\qquad\quad\hspace*{29.2pt}{}\times\exp \biggl\{-\int_0^{\tau} p_t
\bigl(r(t)X_{r(t)^{-2}s}-m_t(\omega) \bigr)\,\d s \biggr\}\dvtx
G_1(\tau) \biggr].
\end{eqnarray*}
Combining the above estimates and using Lemma~\ref{lem9}(i),
we obtain
\[
\gamma_tT(\psi_{\omega}) =e^{-\tau\lambda_1^{\omega}(B_t)+o(1)}
\psi_{\omega}(0)R_0^{\tilde{m}_t(\omega)} \bigl[f(X) \dvtx
X_{[0, t^{\varepsilon_0/2}]}\subset r(t)^{-1}B_t \bigr].
\]
Finally, by the very same argument as for the proof of the
strong localization, it follows that
\[
R_0^{\tilde{m}_t(\omega)} \bigl( X_{[0, t^{\varepsilon_0/2}]}\not\subset
r(t)^{-1}B_t \bigr) =\O \bigl(\exp \bigl\{-(\log
t)^2 \bigr\} \bigr)
\]
on $G_2$. [We have used the second condition in $G_2$ to
control $\psi_{\omega}(0)$.]

Next, we estimate the second term on the
right-hand side of (\ref{decomposition}).
This is rather easy since by the same argument as that for (\ref{F_1-F_2}),
it follows
\begin{eqnarray*}
\bigl|T(\psi_2)\bigr| &\le& \|f\|_{\infty} \exp \bigl\{-(\tau-1) \bigl(
\lambda_1^{\omega}(B_t) +c_3t^{-({\alpha-d+2})/({2\alpha})}
\bigr) \bigr\}
\\
&\le& \|f\|_{\infty} \exp \biggl\{-\tau\lambda_1^{\omega}(B_t)
-\frac{c_3}{2}t^{\varepsilon
_0} \biggr\}.
\end{eqnarray*}
\upqed\end{pf}
Substituting (\ref{lemma10}) into (\ref{replaced}) and
dropping $o(Z_t)$ term, we obtain
\begin{eqnarray*}
&& Q_t \bigl[f(\Tilde{X})\dvtx G_1(t)\cap
G_2 \bigr]
\\
&&\qquad \sim\frac{1}{Z_t} \E \bigl[e^{-t\lambda_1^{\omega}(B_t)} \langle
\phi_{\omega}, 1 \rangle \bigl( \psi_{\omega}(0)R_0^{\tilde{m}_t(\omega)}
\bigl[f(X) \bigr]
\\
&&\hspace*{95.5pt}\qquad\quad{}+ \O \bigl(\exp \bigl\{-(\log t)^2 \bigr\} \bigr)
\bigr)\dvtx G_2 \bigr].
\end{eqnarray*}
The term $\O(\exp\{-(\log t)^2 \})$ is negligible
in view of Lemma~\ref{lem5} and the fact that
$|\langle\phi_{\omega}, 1 \rangle|
\le\|\phi_{\omega}\|_{L^2(B_t)}\|1\|_{L^2(B_t)}$ grows at most
polynomially fast. Therefore we arrive at the expression
%
%
\begin{eqnarray}
\label{almost-done} && Q_t \bigl[f(\Tilde{X})\dvtx
G_1(t)\cap G_2 \bigr]
\nonumber
\\
&&\qquad \sim\int\frac{1}{Z_t} \E \bigl[e^{-t\lambda_1^{\omega}(B_t)} \langle
\phi_{\omega}, 1 \rangle\dvtx G_2, \tilde{m}_t(
\omega) \in\d m \bigr]
\\
&&\qquad\quad{}\times \biggl(\frac{\sqrt{C(d,\alpha)}} {
\sqrt{2}\pi r(t)} \biggr)^{{d}/{4}} \exp \Biggl\{- \sqrt{
\frac{C(d,\alpha)}{2}}|m|^2 \Biggr\} R_0^{m}
\bigl[f(X) \bigr].\nonumber
\end{eqnarray}
Since
$\E[e^{-t\lambda_1^{\omega}(B_t)}
\langle\phi_{\omega}, 1 \rangle
\dvtx G_2, \tilde{m}_t(\omega) \in\d m]$
defines a translation invariant measure on $B(0, (\log t)^{3/4})$,
it is a constant multiple of the Lebesgue measure.
We can determine the constant asymptotically by setting $f=1$,
and it follows that the right-hand side of (\ref{almost-done})
converges to
\[
\int\d m \biggl(\frac{\sqrt{C(d,\alpha)}}{\sqrt{2}\pi} \biggr)^{
{d}/{2}} \exp \Biggl\{-\sqrt{
\frac{C(d,\alpha)}{2}}|m|^2 \Biggr\} R^m_0
\bigl[f(X) \bigr].
\]
\upqed\end{pf*}

\section*{Acknowledgments}
The author would like to thank the referee for a careful reading of the
manuscript and for Remark~\ref{1.1}.



\printaddresses

\end{document}